\documentclass[12pt]{article}
\usepackage{full page,
latexsym, amssymb, amscd, amsthm, graphicx, 
}

 \title{{\bf Logarithmic intertwining operators and associative
algebras}}
 \author{Yi-Zhi Huang and Jinwei Yang}
    \date{}
    \begin{document}
    \bibliographystyle{alpha}
    \maketitle
\newtheorem{thm}{Theorem}[section]
\newtheorem{defn}[thm]{Definition}
\newtheorem{prop}[thm]{Proposition}
\newtheorem{cor}[thm]{Corollary}
\newtheorem{lemma}[thm]{Lemma}
\newtheorem{rema}[thm]{Remark}
\newtheorem{app}[thm]{Application}
\newtheorem{prob}[thm]{Problem}
\newtheorem{conv}[thm]{Convention}
\newtheorem{conj}[thm]{Conjecture}
\newtheorem{cond}[thm]{Condition}
    \newtheorem{exam}[thm]{Example}
\newtheorem{assum}[thm]{Assumption}
     \newtheorem{nota}[thm]{Notation}
\newcommand{\halmos}{\rule{1ex}{1.4ex}}
\newcommand{\pfbox}{\hspace*{\fill}\mbox{$\halmos$}}

\newcommand{\nn}{\nonumber \\}

 \newcommand{\res}{\mbox{\rm Res}}
 \newcommand{\ord}{\mbox{\rm ord}}
\renewcommand{\hom}{\mbox{\rm Hom}}
\newcommand{\edo}{\mbox{\rm End}\ }
 \newcommand{\pf}{{\it Proof.}\hspace{2ex}}
 \newcommand{\epf}{\hspace*{\fill}\mbox{$\halmos$}}
 \newcommand{\epfv}{\hspace*{\fill}\mbox{$\halmos$}\vspace{1em}}
 \newcommand{\epfe}{\hspace{2em}\halmos}
\newcommand{\nord}{\mbox{\scriptsize ${\circ\atop\circ}$}}
\newcommand{\wt}{\mbox{\rm wt}\ }
\newcommand{\swt}{\mbox{\rm {\scriptsize wt}}\ }
\newcommand{\lwt}{\mbox{\rm wt}^{L}\;}
\newcommand{\rwt}{\mbox{\rm wt}^{R}\;}
\newcommand{\slwt}{\mbox{\rm {\scriptsize wt}}^{L}\,}
\newcommand{\srwt}{\mbox{\rm {\scriptsize wt}}^{R}\,}
\newcommand{\clr}{\mbox{\rm clr}\ }
\newcommand{\tr}{\mbox{\rm Tr}}
\newcommand{\C}{\mathbb{C}}
\newcommand{\Z}{\mathbb{Z}}
\newcommand{\R}{\mathbb{R}}
\newcommand{\Q}{\mathbb{Q}}
\newcommand{\N}{\mathbb{N}}
\newcommand{\CN}{\mathcal{N}}
\newcommand{\F}{\mathcal{F}}
\newcommand{\I}{\mathcal{I}}
\newcommand{\V}{\mathcal{V}}
\newcommand{\one}{\mathbf{1}}
\newcommand{\BY}{\mathbb{Y}}
\newcommand{\ds}{\displaystyle}

        \newcommand{\ba}{\begin{array}}
        \newcommand{\ea}{\end{array}}
        \newcommand{\be}{\begin{equation}}
        \newcommand{\ee}{\end{equation}}
        \newcommand{\bea}{\begin{eqnarray}}
        \newcommand{\eea}{\end{eqnarray}}
         \newcommand{\lbar}{\bigg\vert}
        \newcommand{\p}{\partial}
        \newcommand{\dps}{\displaystyle}
        \newcommand{\bra}{\langle}
        \newcommand{\ket}{\rangle}

        \newcommand{\ob}{{\rm ob}\,}
        \renewcommand{\hom}{{\rm Hom}}

\newcommand{\A}{\mathcal{A}}
\newcommand{\Y}{\mathcal{Y}}

\newcommand{\dlt}[3]{#1 ^{-1}\delta \bigg( \frac{#2 #3 }{#1 }\bigg) }

\newcommand{\dlti}[3]{#1 \delta \bigg( \frac{#2 #3 }{#1 ^{-1}}\bigg) }

 \makeatletter
\newlength{\@pxlwd} \newlength{\@rulewd} \newlength{\@pxlht}
\catcode`.=\active \catcode`B=\active \catcode`:=\active
\catcode`|=\active
\def\sprite#1(#2,#3)[#4,#5]{
   \edef\@sprbox{\expandafter\@cdr\string#1\@nil @box}
   \expandafter\newsavebox\csname\@sprbox\endcsname
   \edef#1{\expandafter\usebox\csname\@sprbox\endcsname}
   \expandafter\setbox\csname\@sprbox\endcsname =\hbox\bgroup
   \vbox\bgroup
  \catcode`.=\active\catcode`B=\active\catcode`:=\active\catcode`|=\active
      \@pxlwd=#4 \divide\@pxlwd by #3 \@rulewd=\@pxlwd
      \@pxlht=#5 \divide\@pxlht by #2
      \def .{\hskip \@pxlwd \ignorespaces}
      \def B{\@ifnextchar B{\advance\@rulewd by \@pxlwd}{\vrule
         height \@pxlht width \@rulewd depth 0 pt \@rulewd=\@pxlwd}}
      \def :{\hbox\bgroup\vrule height \@pxlht width 0pt depth
0pt\ignorespaces}
      \def |{\vrule height \@pxlht width 0pt depth 0pt\egroup
         \prevdepth= -1000 pt}
   }
\def\endsprite{\egroup\egroup}
\catcode`.=12 \catcode`B=11 \catcode`:=12 \catcode`|=12\relax
\makeatother

\def\hboxtr{\FormOfHboxtr} 
\sprite{\FormOfHboxtr}(25,25)[0.5 em, 1.2 ex] 

:BBBBBBBBBBBBBBBBBBBBBBBBB | :BB......................B |
:B.B.....................B | :B..B....................B |
:B...B...................B | :B....B..................B |
:B.....B.................B | :B......B................B |
:B.......B...............B | :B........B..............B |
:B.........B.............B | :B..........B............B |
:B...........B...........B | :B............B..........B |
:B.............B.........B | :B..............B........B |
:B...............B.......B | :B................B......B |
:B.................B.....B | :B..................B....B |
:B...................B...B | :B....................B..B |
:B.....................B.B | :B......................BB |
:BBBBBBBBBBBBBBBBBBBBBBBBB |

\endsprite
\def\shboxtr{\FormOfShboxtr} 
\sprite{\FormOfShboxtr}(25,25)[0.3 em, 0.72 ex] 

:BBBBBBBBBBBBBBBBBBBBBBBBB | :BB......................B |
:B.B.....................B | :B..B....................B |
:B...B...................B | :B....B..................B |
:B.....B.................B | :B......B................B |
:B.......B...............B | :B........B..............B |
:B.........B.............B | :B..........B............B |
:B...........B...........B | :B............B..........B |
:B.............B.........B | :B..............B........B |
:B...............B.......B | :B................B......B |
:B.................B.....B | :B..................B....B |
:B...................B...B | :B....................B..B |
:B.....................B.B | :B......................BB |
:BBBBBBBBBBBBBBBBBBBBBBBBB |

\endsprite

\vspace{2em}



\renewcommand{\theequation}{\thesection.\arabic{equation}}
\renewcommand{\thethm}{\thesection.\arabic{thm}}
\setcounter{equation}{0} \setcounter{thm}{0} 
\date{}
\maketitle

\begin{abstract}
We establish an isomorphism between the space of logarithmic 
intertwining operators  among  suitable generalized 
modules for a vertex operator algebra and the space of homomorphisms
between suitable modules for a generalization of Zhu's algebra 
given by Dong-Li-Mason. 
\end{abstract}

\section{Introduction}

In the representation theory of reductive vertex operator algebras 
(vertex operator algebras for which a suitable category of weak
modules is semisimple) and in the construction 
of rational conformal field theories, 
intertwining operators introduced in \cite{FHL}
are in fact the fundamental mathematical objects from which these theories
are developed and constructed.  In \cite{FZ}, for a reductive vertex 
operator algebra 
$V$, Frenkel and Zhu identified 
the spaces of intertwining operators among irreducible $V$-modules 
with suitable spaces constructed 
from (right, bi-, left) modules for Zhu's algebra $A(V)$
associated to the irreducible $V$-modules. See \cite{L} 
for a generalization and a proof of this result. 
This result is very useful for the calculation of fusion rules and 
for the construction of intertwining operators.

To develop the representation theory of vertex operator algebras that
are not reductive, it is necessary to consider certain generalized modules that are 
not completely reducible and the logarithmic intertwining operators among them. 
The theory of logarithmic intertwining operators corresponds to genus-zero 
logarithmic conformal field theories in physics. In fact, 
logarithmic structure in conformal field theory was first
introduced by physicists to describe disorder phenomena \cite{G} and 
logarithmic conformal field theories
have been developed rapidly in recent years. See \cite{HLZ} 
for an introduction and for 
references to the study of logarithmic intertwining operators,
a logarithmic tensor category theory and their connection with 
various works of mathematicians and physicists on logarithmic conformal 
field theories. 

In this general setting, we can ask the following natural question: 
In the case that the generalized
modules involved are
not necessarily completely reducible, can we identify the spaces of 
logarithmic intertwining operators among suitable
generalized modules with some spaces constructed from modules for
certain associative algebras associated to the vertex operator algebra?
We answer this question in the present paper. 
Our answer needs the 
generalizations of Zhu's algebra given by Dong, Li and Mason in \cite{DLM}. 
For a  generalized 
module for the vertex operator algebra, we 
introduce a bimodule for such an associative algebra. 
This bimodule generalizes the bimodule for Zhu's algebra 
given in \cite{FZ}.
Our main result establishes an isomorphism between the space of logarithmic 
intertwining operators  among suitable generalized 
modules and the space of homomorphisms
between suitable modules for a generalization of Zhu's algebra 
given in \cite{DLM}. 
See Theorem \ref{main} for the precise statement of our 
main result. Our method follows the one used in \cite{H2}
and is different from 
the one used in \cite{L}.

Our result will be used in a forthcoming paper on generalized twisted modules 
associated to a not-necessarily-finite-order isomorphism of a vertex operator algebra
(see \cite{H} for the definition and examples of such generalized twisted modules). 
In fact, the results on generalized twisted modules in that forthcoming paper is the main 
motivation for the main theorem that we obtain in this paper. 

The present paper is organized as follows: In the next section, we recall basic 
notions and results on generalized modules for a vertex operator algebra. 
In Section 3, we recall the generalizations of Zhu's algebra by Dong, Li and 
Mason in \cite{DLM}. In Section 4, we introduce and study a bimodule structure 
for such an algebra 
on a quotient of a lower-bounded generalized module for a vertex operator 
algebra.  In Section 5, we begin our study of the relation between 
logarithmic intertwining operators and homomorphisms
between suitable modules for a generalization of Zhu's algebra.
Our main result is stated and proved in Section 6.

\paragraph{Acknowledgments}
The authors are grateful to Haisheng Li for discussions on some calculations
in \cite{DLM}. This research is supported in part by NSF
grant PHY-0901237.

\section{Generalized modules for a vertex operator algebra}

In this paper, we shall assume that the reader is familiar with the basic
notions and results in the theory of vertex operator algebras.
In particular, we assume that the reader is
familiar with weak modules,
$\N$-gradable weak modules, contragredient modules and
related results. Our terminology and conventions follow
those in \cite{FLM}, \cite{FHL} and \cite{H3}.
We shall use $\Z$, $\Z_{+}$, $\N$,  $\R$ and
$\C$ to denote the (sets of) integers, positive integers,
nonnegative integers, real numbers and complex numbers, 
respectively. For $n\in \C$, we use
$\Re{(n)}$ and $\Im{(n)}$ to denote the real and imaginary parts of
$n$.

In this section, we recall the notion of generalized module
for a vertex operator algebra and
related notions in \cite{HLZ} and also some related notions and 
basic properties in \cite{H3}.

We fix a vertex operator algebra $V$ in this paper. 
(In fact, the results in the present paper 
are true for a grading-restricted M\"{o}bius vertex algebra (see \cite{HLZ}).)
We first recall the definition of
generalized $V$-module and related notions in \cite{HLZ} (see also \cite{M}):

\begin{defn}
{\rm A $\C$-graded
vector space $W=\coprod_{n\in \C}W_{[n]}$ equipped with
a linear map
\begin{eqnarray*}
Y_{W}: V\otimes W&\to &W((x))\\
v\otimes w&\mapsto & Y_{W}(v, x)w
\end{eqnarray*}
is called a {\it generalized $V$-module} if it
satisfies all the axioms for $V$-modules except
that  $W$ does not have to satisfy the two grading-restriction
conditions and that instead of requiring $L(0)w=nw$ for $w\in W_{[n]}$, 
we require that the following weaker version of the $L(0)$-grading property,
still called the {\it $L(0)$-grading property}:
For $n\in \C$, the homogeneous
subspaces $W_{[n]}$ are the generalized eigenspaces of $L(0)$ with
eigenvalues $n$, that is, for $n\in \C$, $w\in W_{[n]}$,
there exists $K\in \Z_{+}$, depending on $w$, such that
$(L(0)-n)^{K}w=0$.
We define {\it homomorphisms} (or {\it module maps})
and {\it isomorphisms} (or {\it equivalence})
between generalized $V$-modules, {\it generalized $V$-submodules}
and {\it quotient generalized $V$-modules}
in the obvious way.}
\end{defn}

\begin{defn}
{\rm An {\it irreducible generalized $V$-module} is a generalized $V$-module $W$
such that there is no generalized $V$-submodule of $W$ that is neither
$0$ nor $W$ itself.  A {\it lower bounded generalized $V$-module} is 
a generalized $V$-module $W$ such that $W_{[n]}=0$ when $\Re{(n)}$ is sufficiently negative.
We say that lower-bounded generalized $V$-module $W$ 
{\it has a lowest conformal weight}, or for simplicity, $W$ {\it has a lowest weight}
if there exists $n_{0}\in \C$ such that
$W_{[n_{0}]}\ne 0$ but $W_{[n]}=0$ when $\Re{(n)}<\Re{(n_{0})}$
or $\Re{(n)}=\Re{(n_{0})}$ but $\Im{(n)}\ne \Im{(n_{0})}$. In this case,
we call $n_{0}$, $W_{[n_{0}]}$ and elements of $W_{[n_{0}]}$
the {\it lowest conformal weight} or {\it lowest
weight} of $W$, the {\it lowest weight space}
or {\it lowest weight space}
of $W$ and {\it lowest
conformal weight
vectors} or {\it lowest weight vectors} of $W$, respectively.
A {\it grading restricted generalized $V$-module}  is a generalized $V$-module $W$
such that $W$
is lower bounded
and $\dim W_{[n]}<\infty$
for $n\in \C$.
An {\it (ordinary) $V$-module} is a generalized $V$-module $W$ 
such that $W$ is grading restricted and
$W_{[n]}=W_{(n)}$ for $n\in \C$, where for $n\in \C$, $W_{(n)}$
are the eigenspaces of $L(0)$ with eigenvalues $n$. An {\it irreducible $V$-module}
is a $V$-module such that 
it is irreducible as a generalized $V$-module. 
A {\it generalized $V$-module of length $l$}
is a generalized $V$-module $W$ such that there exist
generalized $V$-submodules $W=W_{1}\supset \cdots
\supset W_{l+1}=0$ such that $W_{i}/W_{i+1}$ for $i=1, \dots, l$
are irreducible $V$-modules.
A {\it finite length generalized $V$-module} is a generalized $V$-module
of length  $l$ for some $l\in \Z_{+}$.
{\it Homomorphisms} and {\it isomorphisms}
between lower-bounded, 
grading-restricted or finite length generalized $V$-modules are
homomorphisms and isomorphisms between the underlying
generalized $V$-modules.}
\end{defn}

\begin{rema}
{\rm By definition, every $V$-module is a
generalized $V$-module, but the converse is not true in general. In particular, 
an irreducible generalized $V$-module in general might not be a $V$-module.
A generalized $V$-module $W$ has a lowest weight if $W$ 
is $\R$-graded and lower-bounded
or if $W$ is lower-bounded and generated by one homogeneous element. 
In particular, the vertex operator algebra $V$, any $\R$-graded 
$V$-module, 
or any irreducible
lower-bounded generalized $V$-module has a lowest weight.}
\end{rema}


Let $W=\coprod_{n\in{\mathbb C}}W_{[n]}$ equipped with $Y_{W}$
be a 
generalized $V$-module. As in \cite{HLZ}, the {\it
opposite vertex operator on $W$
associated to $v\in V$} is defined by
\begin{equation}\label{yo}
Y^o_W(v,x)=Y_W(e^{xL(1)}(-x^{-2})^{L(0)}v,x^{-1}).
\end{equation}

Let $W'$ be
the $\mathbb{C}$-graded vector subspace of $W^{*}$ given by
\begin{equation}\label{Wprime}
W' =  \coprod_{n\in {\mathbb C}} (W_{[n]})^{*}.
\end{equation}
We shall use the notation $\langle\cdot,\cdot\rangle_W$,
or $\langle\cdot,\cdot\rangle$ if the underlying space is clear, to denote
the canonical pairing between $W'$ and $W$.  
As in Section 5.2 of \cite{FHL}, we define a vertex operator map 
$Y'_{W}$ for
$W'$ by
\begin{equation}\label{y'}
\langle Y_{W}'(v,x)w',w\rangle = \langle w', Y^o_W(v,x)w\rangle
\end{equation}
for $v\in V$, $w'\in W'$ and $w\in W$. The correspondence given by $v\mapsto
Y_{W}'(v,x)$ is a linear map {}from $V$ to $({\rm
End}\,W')[[x,x^{-1}]]$. Write
\[
Y_{W}'(v,x)=\sum_{n\in {\mathbb Z}} (Y'_{W})_{n}(v) x^{-n-1}
\]
($(Y'_{W})_{n}(v)\in {\rm End}\,W')$ and 
\[
Y_{W}^{o}(v,x)=\sum_{n\in {\mathbb Z}} (Y^{o}_{W})_{n}(v) x^{-n-1}
\]
($(Y^{o}_{W})_{n}(v)\in {\rm End}\,W)$. Then for $v\in V$, $w'\in W'$ and $w\in W$,

\begin{equation}\label{v'vo}
\langle (Y'_{W})_{n}(v) w', w\rangle = \langle w', (Y^{o}_{W})_{n}(v) w\rangle.
\end{equation}

We have:

\begin{thm}
Let $W$ be a lower-bounded generalized $V$-module.
Then $W'$ equipped with $Y'_{W}$ is a lower-bounded generalized $V$-module.
Moreover, $Y''_{W}|_{V\otimes W}=Y_{W}$.\epf
\end{thm}

The proof of this theorem is the same as those of Theorems 
5.2.1 and 5.3.1 in \cite{FHL}. Note that in this theorem,
$W$ does not have to be grading
restricted. The space $W'$ equipped with $Y_{W}'$ is called the {\it contragredient
of $W$}.

We  also define the
operators $L'(n)$ for $n\in {\mathbb Z}$ by
\[
Y_{W}'(\omega,x)=\sum_{n\in {\mathbb Z}}L'(n)x^{-n-2}.
\]
As in Section 5.2 of \cite{FHL}, by extracting the coefficient of 
$x^{-n-2}$ in (\ref{y'}) with
$v=\omega$ and using the fact that $L(1)\omega=0$, we have
\begin{equation}\label{L'(n)}
\langle L'(n)w',w\rangle=\langle w',L(-n)w\rangle\;\;\mbox{ for
}\;n\in {\mathbb Z}.
\end{equation}

The following fact is useful (see \cite{H3}):

\begin{prop}\label{contrag}
The contragredient of a generalized $V$-module of length $l$
also has length $l$.\epf
\end{prop}

\section{Associative algebras from vertex operator algebras and their modules}

In this section, we recall the generalizations of Zhu's algebra \cite{Z} given
by Dong, Li and Mason in \cite{DLM} associated to a vertex operator algebra.
We prove some elementary but useful results.

Recall that we have fixed our vertex operator algebra $V$
in this paper. 
For $N\in \N$, define a product $*_{N}$ on $V$ by
$$u*_{N}v=\sum_{m=0}^{N}(-1)^{m}{m+N\choose N}\res_{x}
x^{-N-m-1}Y_{V}((1+x)^{L(0)+N}u, x)v$$
for $u, v\in V$. Let $O_{N}(V)$ be the subspace of $V$
spanned by elements of the form
$$\res_{x}x^{-2N-1-n}Y_{V}((1+x)^{L(0)+N}u, x)v$$
for $n\in \Z_{+}$, $u, v\in V$ and of the form
$(L(-1)+L(0))u$ for $u\in V$.

\begin{thm}[\cite{DLM}]
The subspace $O_{N}(V)$ is a two-sided ideal of $V$ under the
product $*_{N}$ and the product $*_{N}$ induces a structure
of associative algebra on the quotient $A_{N}(V)=V/O_{N}(V)$ with the
identity $\one +O_{N}(V)$ and with $\omega +O_{N}(V)$ in the center of
$A_{N}(V)$. \epf
\end{thm}

\begin{rema}
{\rm When $N=0$, $A_{0}(V)$ is the associative algebra first
introduced and
studied by Zhu in  \cite{Z}.}
\end{rema}

Let $W$ be a weak $V$-module and for $N\in \N$, let
\[
\Omega_N(W) = \{w \in W\ |\ (Y_{W})_k(u)w = 0\ {\rm for\ homogeneous}\ u\in V, 
\wt u - k - 1 < -N\}.
\]

\begin{thm}[\cite{DLM}] 
The map $V\to \mbox{\rm End}\; \Omega_{N}(W)$ 
given by $v \mapsto o(v)=(Y_{W})_{{\rm wt}v -
1}(v)$ for homogeneous $v\in V$ induces a structure of $A_N(V)$-module on
$\Omega_N(W)$. \epf
\end{thm}

{}From the commutator formula for vertex operators, 
we know that the space $\hat{V}$ of operators on $V$ of the form 
$(Y_{V})_{n}(u)$ for $u\in V$ and 
$n\in \Z$, equipped with the Lie bracket for operators, is a Lie 
algebra. Giving the operators $(Y_{V})_{n}(u)$
for homogeneous  $u$ the weights $\wt u-n-1$, 
we see that $\hat{V}$ becomes a $\Z$-graded Lie 
algebra. We use $\hat{V}_{(n)}$ to denote the homogeneous subspace of 
weight $n$. Also, $A_{N}(V)$ equipped with the Lie bracket  induced from 
the associative algebra structure is a Lie algebra. 

\begin{prop}[\cite{DLM}]\label{dlm-3}
The map given by $(Y_{V})_{\swt v-1}(v) \mapsto v+O_{N}(V)$ is a surjective 
homomorphism of Lie algebras from $\hat{V}_{(0)}$ to the Lie algebra $A_{N}(V)$.\epf
\end{prop}

Let $W$ be a lower-bounded generalized $V$-module such that 
$W=\coprod_{n \in h_{W}+\N} W_{[n]}$ for some $h_{W}\in \C$
and $W_{[h_{W}]}\ne 0$. 
For $N\in \N$, let 
$$\Omega^{0}_{N}(W)=\coprod_{n=0}^{N} W_{[h_{W}+n]}.$$
It is clear that $\Omega_{N}^{0}(W)\subset \Omega_{N}(W)$. 
Since for $u\in V$, $o(u)$ preserves the weights, $\Omega_{N}^{0}(W)$
is an $A_{N}(V)$-submodule of $\Omega_{N}(W)$.

\begin{rema}\label{congruent}
{\rm A generalized $V$-module $W$ decomposes into generalized submodules
corresponding to the congruence classes of its weights modulo $\Z$.
For $\mu \in \C/\Z$, let
$$W^{\mu} = \coprod_{n \in \mu} W_{[n]}.$$
Then
$$W = \coprod_{\mu \in \C/\Z} W^{\mu}$$
and each $W^{\mu}$ is a generalized $V$-submodule of $W$.  In particular, if a
generalized module $W$ is indecomposable, then 
there exists $h\in \C$ such that $W=\coprod_{n \in h+\Z} W_{[n]}$.
In the case that  $W$ is lower bounded, 
there exists $h_{\mu}\in \C$ for $\mu\in \C/\Z$ such that 
$$W^{\mu} = \coprod_{n \in h_{\mu}+\N} W_{[n]}$$
for $\mu\in \C/\Z$.}
\end{rema}

For a lower-bounded generalized $V$-module $W$, by Remark \ref{congruent},
there exists $h_{\mu}\in \C$ for $\mu\in \C/\Z$ such that 
$$W = \coprod_{\mu \in \C/\Z} W^{\mu}$$
where 
$$W^{\mu} = \coprod_{n \in h_{\mu}+\N} W_{[n]}$$
for $\mu\in \C/\Z$ are lower-bounded generalized $V$-submodules 
of $W$. 
Let 
$$\Omega_{N}^{0}(W)=\sum_{\mu\in \C/\Z}\Omega_{N}^{0}(W^{\mu})\subset W.$$
Since $\Omega_{N}^{0}(W^{\mu})$ is an $A_{N}(V)$-submodule of $\Omega_{N}(W^{\mu})$ 
for each $\mu\in \C/\Z$, $\Omega_{N}^{0}(W)$ is an $A_{N}(V)$-submodule of 
$\Omega_{N}(W)$.

\begin{prop}\label{p2}
Let $W$ be a lower-bounded generalized $V$-module generated by 
$\Omega^{0}_N(W)$ for
some $N\in \N$. 
Then $W$ is spanned by
elements of the form 
$$(Y_{W})_{m_1}(u^1)\cdots (Y_{W})_{m_k}(u^k)w,$$
where $u^1, \dots, u^k$ 
are homogeneous elements of $V$, 
$m_1, \dots, m_k$ are integers such that $\wt  u^i - {m_i} - 1 > 0$  and  $w \in
\Omega^{0}_N(W)$.
\end{prop}
\pf
We know that $W$ is spanned by
elements of the form 
$$(Y_{W})_{m_1}(u^1)\cdots (Y_{W})_{m_k}(u^k)w,$$
where $u^1, \dots, u^k$
are homogeneous elements of $V$, 
$m_1, \dots, m_k$ are integers and  $w \in
\Omega^{0}_N(W)$. We have to show that these elements can be written as 
linear combinations of elements of the same form such that $\wt  u^i - {m_i} - 1 > 0$.
If there exists $i$ such that $\wt u^{i}-m_{i}-1\le 0$, then we can find 
an $i$ such that $\wt u^{i}-m_{i}-1\le 0$ and $\wt u^{j}-m_{j}-1> 0$ for $j>i$. 
The component form of the commutator formula for vertex operators gives
\begin{eqnarray*}
\lefteqn{(Y_{W})_{m_{i}}(u^{i})(Y_{W})_{m_{i+1}}(u^{i+1})
-(Y_{W})_{m_{i+1}}(u^{i+1})(Y_{W})_{m_{i}}(u^{i})}\nn
&&\quad
=\sum_{j\in \N}{m_{i}\choose j}(Y_{W})_{m_{i}+m_{i+1}-j}((Y_{V})_{j}(u^{i})u^{i+1}).
\end{eqnarray*}
Thus we have 
\begin{eqnarray*}
\lefteqn{(Y_{W})_{m_1}(u^1)\cdots (Y_{W})_{m_k}(u^k)w}\nn
&&=
(Y_{W})_{m_1}(u^1)\cdots (Y_{W})_{m_{i-1}}(u^{i-1})\cdot\nn
&&\quad\quad\quad\quad\quad\quad
\cdot (Y_{W})_{m_{i+1}}(u^{i+1})
(Y_{W})_{m_{i}}(u^{i})
(Y_{W})_{m_{i+2}}(u^{i+2})\cdots (Y_{W})_{m_k}(u^k)w\nn
&&\quad +\sum_{j\in \N}{m_{i}\choose j}
(Y_{W})_{m_1}(u^1)\cdots (Y_{W})_{m_{i-1}}(u^{i-1})\cdot\nn
&&\quad\quad\quad\quad\quad\quad
\cdot (Y_{W})_{m_{i}+m_{i+1}-j}((Y_{V})_{j}(u^{i})u^{i+1})
(Y_{W})_{m_{i+2}}(u^{i+2})\cdots (Y_{W})_{m_k}(u^k)w.
\end{eqnarray*}
Using this formula, the fact that $(Y_{W})_{m}(u)\tilde{w}\in 
\Omega^{0}_N(W)$ for homogeneous $u\in V$, $m\in \Z$ and $\tilde{w}\in \Omega^{0}_N(W)$ 
such that $\wt u-m-1\le 0$, and  
inductions on $k$ and on the largest number $i$ such that 
$\wt u^{i}-m_{i}-1\le 0$, we see that 
$$(Y_{W})_{m_1}(u^1)\cdots (Y_{W})_{m_k}(u^k)w,$$
can indeed be written as a
linear combination of elements of the same form such that $\wt  u^i - {m_i} - 1 > 0$.
\epfv

Using the definition of the opposite vertex operators $Y^{o}(u, x)$
for $u\in V$, we see that Proposition \ref{p2} gives:
\begin{cor}\label{p2-cor}
Let $W$ be a lower-bounded generalized $V$-module generated by 
$\Omega^{0}_N(W)$ for
some $N\in \N$.
Then $W$ is spanned by
elements of the form 
$$(Y_{W}^{o})_{m_1}(u^1)\cdots (Y_{W}^{o})_{m_k}(u^k)w,$$
where $u^1, \dots, u^k$ 
are homogeneous elements of $V$, 
$m_1, \dots, m_k$ are integers such that $\wt  u^i - {m_i} - 1 < 0$ and $w \in
\Omega^{0}_N(W)$.\epf
\end{cor}

In the results above, $W$ must be generated by $\Omega^{0}_N(W)$ for
some $N\in \N$. We now show that 
generalized $V$-modules of finite length is lower bounded and 
have this property.

Let $W$ be a generalized $V$-module of length $l$ and $W=W_{1}\supset \cdots\supset 
W_{l+1}=0$ a finite composition series of $W$. Since $W_{i}/W_{i+1}$ for $i=1, \dots, l$ are
irreducible $V$-modules and irreducible $V$-modules as modules 
are lower bounded, there exist homogeneous elements $w_{i}\in W_{i}$ of
weights $h_{i}\in \C$
for $i=1, \dots, l$ such that $w_{i} +W_{i+1}$ for $i=1, \dots, l$ are lowest weight 
vectors of $W_{i}/W_{i+1}$.

\begin{prop}\label{p1}
Let $W$ be a generalized $V$-module of length $l$, $W=W_{1}\supset \cdots\supset 
W_{l+1}=0$ a finite composition series of $W$ and 
$w_{i}\in W_{i}$ homogeneous elements of
weights $h_{i}\in \C$
for $i=1, \dots, l$ such that $w_{i} +W_{i+1}$ for $i=1, \dots, l$ are lowest weight 
vectors of $W_{i}/W_{i+1}$. Let $N$ be a positive integer such
that $|\Re(h_i) - \Re(h_j)| \leq N$ for  $i \neq j$, $i, j\in \{1,
\dots, l\}$ and $r=\min_{i\in \{1, \dots, l\}}\Re(h_i)$. 
Then $W$ is lower bounded,
the real number $r$ is the smallest real part of the weights of elements
of $W$ and the subset $\{w_{1}, \dots, w_{l}\}$ of $W$ is in $\Omega^{0}_{N}(W)$
and generates $W$. In particular, $\Omega^{0}_{N}(W)$ generates
$W$.
\end{prop}
\pf
Since $w_{i} +W_{i+1}$ is a lowest weight 
vector of the irreducible $V$-modules $W_{i}/W_{i+1}$ for $i=1, \dots, l$,
$W$ as a graded vector space is isomorphic to 
$\coprod_{i=1}^{l}W_{i}/W_{i+1}$. 
Since the lowest weight of $W_{i}/W_{i+1}$
is $h_{i}$ for $i=1, \dots, l$, the real part of the weight of any homogeneous 
vector of $W$ is larger than or equal to $r=\min_{i\in \{1, \dots, l\}}\Re(h_i)$. 
So $W$ is lower bounded,
$r$ is the smallest real part of the  weights of the elements of 
the graded space $\coprod_{i=1}^{l}W_{i}/W_{i+1}$ and thus
$r$ is also the smallest real part of the weights of the elements of $W$.

Since $|\Re(h_i) - \Re(h_j)| \leq N$, $w_{1}, \dots, w_{l}\in \coprod_{\Re(n)\le r+N}
W_{[n]}$. By definition, we know that $\coprod_{\Re(n)\le r+N}
W_{[n]}\subset \Omega_{N}^{0}(W)$. Thus $w_{1}, \dots, w_{l}\in \Omega_{N}^{0}(W)$.

Let $\widetilde{W}$ be the
generalized $V$-submodule generated by $w_{i}$ for $i=1, \dots, l$.
Since
$W_{i}/W_{i+1}$ for $i=1, \dots, l$ are irreducible,
$w_{i}+W_{i+1}$ for $i=1, \dots, l$ are generators of
$W_{i}/W_{i+1}$.  
We now show that $W=\widetilde{W}$. Since $W_{l}=W_{l}/W_{l+1}$ is
generated by $w_{l}$, we see that $W_{l}\subset \widetilde{W}$.
Now assume that $W_{m}\subset \widetilde{W}$. Then
since $W_{m-1}/W_{m}$ is generated by $w_{m-1}+W_{m}$,
every element of $W_{m-1}/W_{m}$ is a linear combination of
elements of the form
\begin{eqnarray*}
\lefteqn{(Y_{W_{m-1}/W_{m}})_{n_{1}}(u^{1})\cdots (Y_{W_{m-1}/W_{m}})_{n_{k}}(u^{k})
(w_{m-1}+W_{m})}\nn
&&=(Y_{W_{m-1}})_{n_{1}}(u^{1})\cdots (Y_{W_{m-1}})_{n_{k}}(u^{k})w_{m-1}+W_{m}.
\end{eqnarray*}
Thus elements of $W_{m-1}$ are linear combinations of elements of
the form 
$$(Y_{W_{m-1}})_{n_{1}}(u^{1})\cdots (Y_{W_{m-1}})_{n_{k}}(u^{k})w_{m-1}+w$$
where
$w\in W_{m}$. Since 
$$(Y_{W_{m-1}})_{n_{1}}(u^{1})\cdots (Y_{W_{m-1}})_{n_{k}}(u^{k})
w_{m-1}\in \widetilde{W}$$
and $w\in W_{m}\subset \widetilde{W}$,
$$(Y_{W_{m-1}})_{n_{1}}(u^{1})\cdots (Y_{W_{m-1}})_{n_{k}}(u^{k})w_{m-1}+w\in \widetilde{W}.$$
So $W_{m-1}\subset \widetilde{W}$. By the principle of induction,
$W=W_{1}\subset \widetilde{W}$. Thus we see that $\{w_{1}, \dots, w_{l}\}$
generates $W$.
\epfv

\begin{rema}
{\rm Let $W$ be a generalized $V$-module of finite length. 
By Propositions \ref{contrag} and \ref{p1}, the contragredient module 
$W'$ is generated 
by $\Omega^{0}_{N}(W')$ for some $N\in \N$. }
\end{rema}


\section{$A_{N}(V)$-bimodules from generalized $V$-modules}

In this section, for a generalized $V$-module $W$ and $N\in \N$,
we introduce an $A_{N}(V)$-bimodule $A_{N}(W)$.
These bimodules should be viewed as generalizations of bimodules for Zhu's
algebra introduced in \cite{FZ}. We emphasize that the formulas defining 
the right actions on the bimodules 
given in this section are different from but equivalent to 
the one in \cite{FZ} in the case
$N=0$. Our formulas are more natural and conceptual.

In this section, we fix a generalized $V$-module $W$. 
We need the semisimple part $L(0)_s \in {\rm End}\; W$ of the operator $L(0)$
on $W$ defined by
$$L(0)_sw  =  nw$$
for $w\in W_{[n]}$, $n \in \mathbb{C}$.
Recall from \cite{HLZ} that 
we have the commutator formula
\begin{eqnarray*}
[L(0)_{s}, Y_{W}(u, x_{0})]&=&[L(0), Y_{W}(u, x_{0})]\nn
&=&Y_{W}(L(0)u, x_{0})
+x_{0}\frac{d}{dx_{0}}Y_{W}(u, x_{0})
\end{eqnarray*}
for $u\in V$. In particular, we have
$$[L(0)_{s}, L(-1)]=L(-1).$$
Thus we have the $L(0)_{s}$-conjugation property
\begin{equation}\label{conj-y-w}
y^{L(0)_{s}}Y_{W}(u, x)y^{-L(0)_{s}}
=Y_{W}(y^{L(0)_{s}}u, xy)
\end{equation}
for $u\in V$ and
\begin{equation}\label{conj-l-1}
y^{L(0)_{s}}e^{xL(-1)}y^{-L(0)_{s}}=e^{xyL(-1)}.
\end{equation}
These formulas for $L(0)_{s}$ certainly hold also for $L(0)$. 

We also need the map 
\begin{eqnarray*}
Y_{WV}^{W}: W \otimes V &\longrightarrow& W[[x, x^{-1}]]\\
w\otimes u &\longmapsto& Y_{WV}^{W}(w, x)u
\end{eqnarray*}
defined  in \cite{FHL} by
$$Y_{WV}^{W}(w, x)u=e^{xL(-1)}Y_{W}(u, -x)w$$
for $u\in V$ and $w\in W$. 
Proposition 5.1.2 and Remark 5.4.2 in \cite{FHL} give in particular the following:

\begin{prop}[\cite{FHL}]\label{t1}
The map $Y_{WV}^W$ is an intertwining 
operator of type ${W\choose WV}$. In particular, the Jacobi identity
\begin{eqnarray}\label{jacobi-1}
&{\displaystyle x_{0}^{-1}\delta\left(\frac{x_{1}-x_{2}}{x_{0}}\right)
Y_W(u,x_{1})Y_{WV}^W(w,x_{2})v -x_{0}^{-1}\delta\left(\frac{x_{2}-x_{1}}{-x_{0}}\right)
Y_{WV}^W(w, x_{2})Y_{V}(u,x_{1})v }&\nn
&{\displaystyle = x_{2}^{-1}\delta\left(\frac{x_{1}-x_{0}}{x_{2}}\right)
Y_{WV}^W(Y_W(u, x_{0})w, x_{2})v}&
\end{eqnarray}
holds for $u, v\in V$ and $w\in W$. Moreover, the Jacobi identity 
\begin{eqnarray}\label{jacobi-2}
&{\displaystyle x_{0}^{-1}\delta\left(\frac{x_{1}-x_{2}}{x_{0}}\right)
Y_{WV}^W(w, x_{1})Y_{V}(v, x_{2})u - 
x_{0}^{-1}\delta\left(\frac{x_{2}-x_{1}}{-x_{0}}\right)Y_W(v, x_{2})Y_{WV}^W(w, x_{1})u} \nn
&{\displaystyle = x_{2}^{-1}\delta\left(\frac{x_{1}-x_{0}}{x_{2}}\right)
Y_{WV}^W(Y_{WV}^W(w, x_{0})v, x_{2})u}&
\end{eqnarray}
also holds for $u, v\in V$ and $w\in W$.
\end{prop}

We also have:

\begin{prop}
For $w\in W$, 
$$y^{L(0)}Y_{WV}^{W}(w, x)y^{-L(0)}
=Y_{WV}^{W}(y^{L(0)}w, xy).$$
This formula also holds with $L(0)$ replaced by $L(0)_{s}$. 
\end{prop}
\pf
These follow from the definition of $Y_{WV}^{W}$ and (\ref{conj-l-1}) for $L(0)$ or
$L(0)_{s}$.
\epfv

For $N\in \N$, $u \in V$ and $w \in
W$, we define
\begin{eqnarray*}
u *_N w &=& \sum_{m = 0}^N (-1)^m {m+N\choose N}
\res_x x^{-N - m - 1}Y_{W}((1+x)^{L(0)+N}u, x)w,\\
w *_N u &=& \sum_{m = 0}^N (-1)^m {m+N\choose N}
\res_x x^{-N - m - 1}\cdot\nn
&&\quad\quad\quad\quad \cdot (1+x)^{-(L_{W}(-1)+L_{W}(0))}
Y_{WV}^{W}((1+x)^{L(0)+N}w, x)u.
\end{eqnarray*}
Let $O_N(W)$ be the subspace of $W$ spanned by elements of the form
$$u \circ_N w=\res_x x^{-2N -2}Y_{W}((1+x)^{L(0)+N}u, x)w$$
for
$u \in V$ and $w \in W$. Let $A_N(W) = W/O_N(W)$.

\begin{rema}
{\rm Note that in the case of $N=0$, our right action is different from the 
right action in \cite{FZ}. Certainly, these right actions induce the same 
right action on $A_0(W)$ (see Remark \ref{alt-r-action} below). 
The advantage of defining the right action 
using the formula above is that many formulas involving the right action
can be proved in the same way as the proofs of the corresponding formulas for the algebra 
or for the left action, with
some of the vertex operator maps replaced by the map $Y_{WV}^{W}$.
On the other hand, we note that
the definition of the right action above makes sense only for generalized 
$V$-modules, not for weak $V$-modules that are not generalized $V$-modules. }
\end{rema}

The following lemma generalizes Lemma 2.1 in \cite{DLM}:

\begin{lemma}\label{l1}
Let $u \in V$ and $w \in W$. 

\begin{enumerate}

\item We have
\begin{eqnarray*}
u \ast_N w &=& \sum_{m = 0}^N {m+N\choose N}(-1)^N
{\rm Res}_x x^{-N-m-1} \cdot\nn
&&\quad\quad\quad\cdot (1+x)^{-(L(-1)+L(0))}
Y_{WV}^W((1 + x)^{L(0)+m-1}w, x)u,\\
w \ast_N u &=& \sum_{m = 0}^N {m+N\choose N}(-1)^N {\rm Res}_x 
x^{-N-m-1} Y_W((1+ x)^{L(0)+m-1}u, x)w  .
\end{eqnarray*}

\item For $p \geq q \geq
0$, we have
\begin{eqnarray*} 
{\rm Res}_x x^{-2N - 2 - p}\:Y_W((1
+x)^{L(0)+N+q}u,x)w &\in& O_N(W),\\
{\rm Res}_x x^{-2N - 2 - p}\:(1+x)^{-(L(-1)+L(0))}
Y_{WV}^W((1 + x)^{L(0)+N+q}w,x)u &\in& O_N(W).
\end{eqnarray*}
In particular, 
$$w \circ_N u=\res_x x^{-2N -2}
(1+x)^{-(L(-1)+L(0))}Y_{WV}^{W}((1+x)^{L(0)+N}w, x)u\in O_N(W).$$

\item We have 
\begin{eqnarray*}
u \ast_N w - w \ast_N u &=& {\rm Res}_xY_W((1 + x)^{L(0)-1}u, x)w,\\
w \ast_N u - u \ast_N w &=&{\rm Res}_x
(1+x)^{-(L(-1)+L(0))}Y_{WV}^W((1 + x)^{L(0)-1}w,
x)u.
\end{eqnarray*}
\end{enumerate}
\end{lemma}
\pf 
Using the definition of $Y_{WV}^{W}$, we obtain
\begin{eqnarray*}
Y_W(u, x)w 
&=& (1+x)^{L(-1)+L(0)}Y_{WV}^W\left((1 + x)^{-L(0)}w,
\frac{-x}{1 +
x}\right)(1 + x)^{-L(0)}u,\\
Y_{WV}^W(w, x)u &=& (1+x)^{L(-1)+L(0)}Y_{W}\left((1 + x)^{-L(0)}u,
\frac{-x}{1 + x}\right)(1 + x)^{-L(0)}w.
\end{eqnarray*}
Then 
\begin{eqnarray*} 
u \ast_N w &=& \sum_{m = 0}^N (-1)^m
{m+N\choose N}
\res_y y^{-N-m-1}Y_{W}((1 + y)^{L(0)+N}u, y)w\nn
&=& \sum_{m = 0}^N (-1)^m {m+N\choose N}\cdot\nn
&&\quad\quad\quad\cdot 
\res_y y^{-N-m-1}(1+y)^{L(-1)+L(0)}
Y_{WV}^W\left((1 + y)^{-L(0)+N}w, \frac{-y}{1 + y}\right)u\nn
&=& \sum_{m = 0}^N (-1)^N {m+N\choose N}
\res_x x^{-N-m-1}(1+x)^{-(L(-1)+L(0))}Y_{WV}^W((1 + x)^{L(0)+m-1}w, x)u
\end{eqnarray*}
and
\begin{eqnarray*} 
w \ast_N u &=& \sum_{m = 0}^N (-1)^m
{m+N\choose N}
\res_y y^{-N-m-1}(1+y)^{-(L(-1)+L(0))}Y_{WV}^W((1 + y)^{L(0)+N}w, y)u\nn
&=& \sum_{m = 0}^N (-1)^m {m+N\choose N}
\res_y y^{-N-m-1}
Y_{W}\left((1 + y)^{-L(0)+N}u, \frac{-y}{1 + y}\right)w
 \nn
&=& \sum_{m = 0}^N (-1)^N {m+N\choose N}
\res_x x^{-N-m-1}Y_{W}((1+ x)^{L(0)+m-1}u, x)w,
\end{eqnarray*}
where in both formulas, the last steps are obtained by changing the 
variable $x = \frac{-y}{1+y}$. This
proves Part 1.

Similarly, 
\begin{eqnarray*} 
w \circ_N u &=& 
\res_y y^{-2N-2}(1+y)^{-(L(-1)+L(0))}Y_{WV}^W((1 + y)^{L(0)+N}w, y)u\nn
&=& 
\res_y y^{-2N-2}
Y_{W}\left((1 + y)^{-L(0)+N}u, \frac{-y}{1 + y}\right)w
\nn
&=& \res_x x^{-2N-2}Y_{W}((1+ x)^{L(0)+N}u, x)w\nn
&\in & O_N(W).
\end{eqnarray*}
Now the proof of Part 2 is similar to the proof of Lemma 2.1.2 of
\cite{Z}.

Using Part 1,  we obtain
\begin{eqnarray*}
\lefteqn{u \ast_N w - w \ast_N u}\nn
&&= {\rm Res}_x \left(\sum_{m = 0}^N {m+N\choose N} \frac{(-1)^m(1 +
x)^{N + 1} - (-1)^N(1 + x)^m}{x^{N + m + 1}}\right)
Y_W((1 + x)^{L(0)_{s} - 1}u, x)w\nn
&&={\rm Res}_x Y_W((1 + x)^{L(0)_{s} - 1}u, x)w,
\end{eqnarray*}
where the last step uses the formula
$$ \sum_{m = 0}^N
{m+N\choose N} \frac{(-1)^m(1 + x)^{N + 1} - (-1)^N(1 + x)^m}{x^{N +
m + 1}} = 1$$
given by Proposition 5.2 in \cite{DLM}.
This proves the first property in Part 3. The second
is similar.
\epfv

\begin{rema}\label{alt-r-action}
{\rm By Part 1 in Lemma \ref{l1}, we see that we can also define 
the right action $*_{N}$ of $V$ on $W$ by 
$$w*_{N}u=\sum_{m = 0}^N {m+N\choose N}(-1)^N {\rm Res}_x 
x^{-N-m-1} Y_W((1+ x)^{L(0)+m-1}u, x)w$$
for $u\in V$ and $w\in W$. The advantage of this right action is that $W$ 
does not have to be 
a generalized $V$-module.}
\end{rema}

\begin{lemma}\label{l2}
The subspace $O_N(W)$ of $W$ is invariant under the left and right actions of $V$
above. 
\end{lemma}
\pf 
We need to prove 
\begin{eqnarray}
&\label{a2}(u \circ_N w) \ast_N v \in O_N(W),\\
&\label{a2.5}v \ast_N (u \circ_N w) \in O_N(W)
\end{eqnarray}
for $u, v\in V$ and $w\in W$.

We prove (\ref{a2}) first. For homogeneous $u\in V$ and $w\in W$, we have
\begin{eqnarray*}
u \circ_N w &=& \res_y y^{-2N - 2}Y_{W}((1+y)^{L(0)+N}u,
y)w\nn
&=& \res_y y^{-2N - 2}(1+y)^{\swt u+N}Y_{W}(u,
y)w.
\end{eqnarray*}
For $L(0)$ acting on $W$, we write $L(0)_{n}=L(0)-L(0)_{s}$.
Then for homogeneous $u\in V$, $v\in V$ and homogeneous $w\in W$, 
\begin{eqnarray*}
\lefteqn{(u \circ_N w) \ast_N v} \nn
&&=\sum_{m = 0}^N(-1)^m{m+N \choose N}{\rm Res}_{x_{2}}{\rm Res}_{y}
y^{-2N -2}x_{2}^{-N-m-1}(1+y)^{\swt u+N}\cdot\nn
&&\quad\quad\quad\quad\quad\quad
 \cdot (1+x_{2})^{-(L(-1)+L(0))}
Y_{WV}^W((1 + x_{2})^{L(0)+N}Y_W(u, y)w, x_{2})v\nn
&&=\sum_{m = 0}^N(-1)^m{m+N \choose N}{\rm Res}_{x_{2}}{\rm Res}_{y}
y^{-2N -2}x_{2}^{-N-m-1}(1+y)^{\swt u+N}(1 + x_{2})^{\swt u+\swt w+N}\cdot\nn
&&\quad\quad\quad\quad\quad\quad\cdot 
(1+x_{2})^{-(L(-1)+L(0))}Y_{WV}^W(Y_W(u, (1+x_{2})y)
(1+x_{2})^{L(0)_{n}}w, x_{2})v\nn
&&=\sum_{m = 0}^N(-1)^m{m+N \choose N}{\rm Res}_{x_{2}}{\rm Res}_{x_{0}}
x_{0}^{-2N -2}x_{2}^{-N-m-1}(1+x_{2}+x_{0})^{\swt u+N}
\cdot\nn
&&\quad\quad\quad\quad\quad\quad\cdot (1 + x_{2})^{\swt w+2N+1}
(1+x_{2})^{-(L(-1)+L(0))}
Y_{WV}^W(Y_W(u, x_{0})(1+x_{2})^{L(0)_{n}}w, x_{2})v\nn
&&=\sum_{m = 0}^N(-1)^m{m+N \choose N}{\rm Res}_{x_{2}}{\rm
Res}_{x_{0}}\res_{x_{1}}x_{0}^{-1}\delta\left(\frac{x_{1}-x_{2}}{x_{0}}\right)
x_{0}^{-2N -2}x_{2}^{-N-m-1}
\cdot\nn
&&\quad\quad\quad\quad\quad\quad\cdot (1 + x_{2})^{\swt w +2N+1}
(1+x_{2}+x_{0})^{\swt u+N}
(1+x_{2})^{-(L(-1)+L(0))}\cdot\nn
&&\quad\quad\quad\quad\quad\quad\cdot Y_{W}(u, x_{1})Y_{WV}^W((1+x_{2})^{L(0)_{n}}w, x_{2})v\nn
&&\quad - \sum_{m = 0}^N(-1)^m{m+N \choose N}{\rm Res}_{x_{2}}{\rm
Res}_{x_{0}}\res_{x_{1}}x_{0}^{-1}\delta\left(\frac{x_{2}-x_{1}}{-x_{0}}\right)
x_{0}^{-2N -2}x_{2}^{-N-m-1}\cdot\nn
&&\quad\quad\quad\quad\quad\quad\cdot
(1+x_{2}+x_{0})^{\swt u+N}(1 + x_{2})^{\swt w +2N+1}
(1+x_{2})^{-(L(-1)+L(0))}
\cdot\nn
&&\quad\quad\quad\quad\quad\quad\cdot Y_{WV}^W((1+x_{2})^{L(0)_{n}}w, x_{2})Y_{V}(u, x_{1})v\nn
&&=\sum_{m = 0}^N(-1)^m{m+N \choose N}{\rm Res}_{x_{2}}{\rm
Res}_{x_{0}}\res_{x_{1}}x_{0}^{-1}\delta\left(\frac{x_{1}-x_{2}}{x_{0}}\right)
(x_{1}-x_{2})^{-2N -2}x_{2}^{-N-m-1}\cdot\nn
&&\quad\quad\quad\quad\quad\quad\cdot
(1+x_{1})^{\swt u+N}(1 + x_{2})^{\swt w +2N+1}
(1+x_{2})^{-(L(-1)+L(0))}\cdot\nn
&&\quad\quad\quad\quad\quad\quad\cdot Y_{W}(u, x_{1})Y_{WV}^W((1+x_{2})^{L(0)_{n}}w, x_{2})v\nn
&&\quad - \sum_{m = 0}^N(-1)^m{m+N \choose N}{\rm Res}_{x_{2}}{\rm
Res}_{x_{0}}\res_{x_{1}}x_{0}^{-1}\delta\left(\frac{x_{2}-x_{1}}{-x_{0}}\right)
(-x_{2}+x_{1})^{-2N -2}x_{2}^{-N-m-1}\cdot\nn
&&\quad\quad\quad\quad\quad\quad\cdot
(1+x_{1})^{\swt u+N}(1 + x_{2})^{\swt w +2N+1}
(1+x_{2})^{-(L(-1)+L(0))}\cdot\nn
&&\quad\quad\quad\quad\quad\quad\cdot
Y_{WV}^W((1+x_{2})^{L(0)_{n}}w, x_{2})Y_{V}(u, x_{1})v\nn
&&=\sum_{m = 0}^N\sum_{i\in \N}(-1)^{m+i}{m+N \choose N}{-2N
- 2 \choose i}{\rm Res}_{x_{2}}
\res_{x_{1}}
x_{1}^{-2N -2-i}x_{2}^{-N-m-1+i}\cdot\nn
&&\quad\quad\quad\quad\quad\quad\cdot (1 + x_{2})^{\swt w +2N+1}
(1+x_{2})^{-(L(-1)+L(0))}\cdot\nn
&&\quad\quad\quad\quad\quad\quad\cdot 
Y_{W}((1+x_{1})^{L(0)+N}u, x_{1})Y_{WV}^W((1+x_{2})^{L(0)_{n}}w, x_{2})v\nn
&&\quad - \sum_{m = 0}^N\sum_{i\in \N}(-1)^{m+i}{m+N \choose N}{-2N
- 2 \choose i}{\rm Res}_{x_{2}}
\res_{x_{1}}
x_{2}^{-3N -m-3-i}x_{1}^{i}(1+x_{1})^{\swt u+N}\cdot\nn
&&\quad\quad\quad\quad\quad\quad\cdot
(1+x_{2})^{-(L(-1)+L(0))}
Y_{WV}^W((1 + x_{2})^{L(0) +2N+1}w, x_{2})Y_{V}(u, x_{1})v,
\end{eqnarray*}
where in the third step, we have changed the variable $y=\frac{x_{0}}{1+x_{2}}$
and in the fourth step, we have used the Jacobi identity (\ref{jacobi-1}).
Since all the terms in the first sum of the right-hand side lie in $O_N(W)$ by definition
and all the terms in the second sum of the right-hand side lie in $O_N(W)$ by 
Part 2 of Lemma \ref{l1}, 
(\ref{a2}) holds.

For homogeneous $u, v\in V$ and $w\in W$,
using Part 3 in Lemma \ref{l1} and (\ref{a2}), we obtain
\begin{eqnarray*}
\lefteqn{v \ast_N (u \circ_N w) }\nn
&&=(u \circ_N w) \ast_N v+\res_{x_{1}}(1 +
x_{1})^{\swt v - 1} Y_W(v, x_{1})(u \circ_N w)\nn
&&\equiv \res_{x_{1}}(1 +
x_{1})^{\swt v - 1} Y_W(v, x_{1})(u \circ_N w)\ \ {\rm mod}\ O_N(W)\nn
&&= \res_{x_{1}}\res_{x_{2}}(1 +
x_{1})^{\swt v - 1}(1 + x_{2})^{\swt u+N}x_{2}^{-2N - 2}
Y_W(v, x_{1})Y_W(u, x_{2})w\nn
&&=\res_{x_{1}}\res_{x_{2}}(1 +
x_{1})^{\swt v - 1}(1 + x_{2})^{\swt u+N}x_{2}^{-2N - 2}
Y_W(u, x_{2})Y_W(v, x_{1})w\nn
&&\quad +\res_{x_{1}}\res_{x_{2}}\res_{x_{0}}
x_{2}^{-1}\delta\left(\frac{x_{1}-x_{0}}{x_{2}}\right)\cdot\nn
&&\quad\quad\quad\quad\quad\quad \cdot(1 +
x_{1})^{\swt v - 1}(1 + x_{2})^{\swt u+N}x_{2}^{-2N - 2}
Y_W(Y_{V}(v, x_{0})u, x_{2})w\nn
&&\equiv \res_{x_{2}}\res_{x_{0}}(1 +
x_{2}+x_{0})^{\swt v - 1}(1 + x_{2})^{\swt u+N}x_{2}^{-2N - 2}
Y_W(Y_{V}(v, x_{0})u, x_{2})w\ \ {\rm mod}\ O_N(W)\nn
&&= \sum_{i \in \N}{{\rm wt}\;v - 1 \choose i}\res_{x_{2}}
(1 + x_{2})^{\swt u+\swt v+N - 1 - i}x_{2}^{-2N - 2}Y_W(Y_i(v)u, x_{2})w\nn
&&=\sum_{i \in \N}{{\rm wt}\;v - 1 \choose i}\res_{x_{2}}
x_{2}^{-2N - 2}Y_W((1 + x_{2})^{L(0)+N}Y_i(v)u, x_{2})w\nn
&&\in O_N(W),
\end{eqnarray*}
proving (\ref{a2.5}).
\epfv

The main result in this section is the following:

\begin{thm}
The left and right actions of $V$ on $W$ induce an
$A_N(V)$-bimodule structure on $A_N(W)$.
\end{thm}
\pf 
Lemma \ref{l2} says that the left and right actions of $V$ on $W$ give
left and right actions of $V$ on $A_{N}(W)$. 
We first need to show that these left and right actions of $V$ on $A_{N}(W)$
in fact give left and right actions of $A_{N}(V)$ on $A_{N}(W)$, that is, 
we need to prove
\begin{eqnarray}
&\label{a5}(L(-1)u + L(0)u) \ast_N w \in O_N(W),\\
&\label{a5.5}w \ast_N (L(-1)u + L(0)u) \in O_N(W),\\
&(u \circ_N v) \ast_N w \in O_N(W),\\
&\label{a6}w \ast_N (v \circ_N u) \in O_N(W)
\end{eqnarray}
for $u, v\in V$ and $w\in W$. The proof of these formulas are similar to
the proof of Lemma \ref{l2} and we omit them. 

Next we need to prove that these left and right actions indeed give
left and right $A_{N}(V)$ modules, that is, we need to prove
\begin{eqnarray}
\label{a3}u \ast_N (v \ast_N w) &\equiv& (u \ast_N v) \ast_N w\ {\rm mod}\  O_N(W),\\
\label{a4}w \ast_N (v \ast_N u) &\equiv& (w \ast_N v) \ast_N u \ {\rm mod}\  O_N(W)
\end{eqnarray}
$u, v\in V$ and $w\in W$.
We prove only (\ref{a4}) here; the proof of (\ref{a3}) is similar.

For $v\in V$ and homogeneous $w\in W$, 
\begin{eqnarray*}
w \ast_N v &=& \sum_{m=0}^N(-1)^m{m+N \choose N}\res_y
y^{-N - m - 1}(1+y)^{-(L(-1)+L(0))}Y_{WV}^W((1+y)^{L(0)+N}w,
y)v\nn
&=& \sum_{m=0}^N(-1)^m{m+N \choose N}\res_y
y^{-N - m - 1}(1+y)^{\swt w+N}\cdot\nn
&&\quad\quad\quad\quad\quad\quad\quad\quad\cdot
(1+y)^{-(L(-1)+L(0))}
Y_{WV}^W((1+y)^{L(0)_{n}}w,
y)v.
\end{eqnarray*}
Then for $u\in V$, homogeneous $v\in V$ and homogeneous $w\in W$, 
\begin{eqnarray*}
\lefteqn{(w\ast_N v)\ast_N u}\nn
&&=\sum_{m = 0}^N\sum_{n = 0}^N(-1)^{m+n}{m+N \choose N}{n+N
\choose N}\res_{x_{2}}\res_{y}
y^{-N-m-1}x_{2}^{-N-n-1}(1+y)^{\swt w + N}\cdot\nn
&&\quad\quad\quad\quad \cdot
(1+x_{2})^{-(L(-1)+L(0))}
Y_{WV}^W((1+x_{2})^{L(0)+ N}(1+y)^{-(L(-1)+L(0))}\cdot\nn
&&\quad\quad\quad\quad\cdot
Y_{WV}^W((1+y)^{L(0)_{n}}w, y)v,
x_{2})u\nn
&&=\sum_{m = 0}^N\sum_{n = 0}^N(-1)^{m+n}{m+N \choose N}{n+N
\choose N}\res_{x_{2}}\res_{y}
y^{-N-m-1}x_{2}^{-N-n-1}(1+y)^{\swt w + N}\cdot\nn
&&\quad\quad\quad\quad \cdot
(1+x_{2})^{-(L(-1)+L(0))}
Y_{WV}^W((1+y)^{-((1+x_{2})L(-1)+L(0))}(1+x_{2})^{L(0)+ N}\cdot\nn
&&\quad\quad\quad\quad\cdot
Y_{WV}^W((1+y)^{L(0)_{n}}w, y)v,
x_{2})u\nn
&&=\sum_{m = 0}^N\sum_{n = 0}^N(-1)^{m+n}{m+N \choose N}{n+N
\choose N}\res_{x_{2}}\res_{y}
y^{-N-m-1}x_{2}^{-N-n-1}(1+y)^{\swt w + N}\cdot\nn
&&\quad\quad\quad\quad \cdot
(1+x_{2})^{-(L(-1)+L(0))}
Y_{WV}^W((1+y)^{-(x_{2}L(-1)+L(0))}e^{-yL(-1)}(1+x_{2})^{L(0)+ N}\cdot\nn
&&\quad\quad\quad\quad\cdot
Y_{WV}^W((1+y)^{L(0)_{n}}w, y)v,
x_{2})u\nn
&&=\sum_{m = 0}^N\sum_{n = 0}^N(-1)^{m+n}{m+N \choose N}{n+N
\choose N}\res_{x_{2}}\res_{y}
y^{-N-m-1}x_{2}^{-N-n-1}(1+y)^{\swt w + N}\cdot\nn
&&\quad\quad\quad\quad \cdot
(1+x_{2})^{-(L(-1)+L(0))}(1+y)^{-(L(-1)+L(0))}
Y_{WV}^W((1+x_{2})^{L(0)+ N}\cdot\nn
&&\quad\quad\quad\quad\cdot Y_{WV}^W((1+y)^{L(0)_{n}}w, y)v,
x_{2})(1+y)^{L(-1)+L(0)}u\nn
&&=\sum_{m = 0}^N\sum_{n = 0}^N(-1)^{m+n}{m+N \choose N}{n+N
\choose N}\res_{x_{2}}\res_{y}
y^{-N-m-1}x_{2}^{-N-n-1}(1+y)^{\swt w + N}\cdot\nn
&&\quad\quad\quad\quad\cdot
(1+x_{2})^{\swt w+\swt v+ N}
(1+x_{2})^{-(L(-1)+L(0))}(1+y)^{-(L(-1)+L(0))}
\cdot\nn
&&\quad\quad\quad\quad\cdot
Y_{WV}^W(Y_{WV}^W((1+x_{2}+(1+x_{2})y)^{L(0)_{n}}w, (1+x_{2})y)v,
x_{2})u\nn
&&=\sum_{m = 0}^N\sum_{n = 0}^N(-1)^{m+n}{m+N \choose N}{n+N
\choose N}\res_{x_{2}}\res_{x_{0}}
x_{0}^{-N-m-1}x_{2}^{-N-n-1}(1+x_{2})^{\swt v+ N+m}\cdot\nn
&&\quad\quad\quad\quad\cdot
(1+x_{2}+x_{0})^{\swt w + N} (1+x_{2}+x_{0})^{-(L(-1)+L(0))}
\cdot\nn
&&\quad\quad\quad\quad\cdot
Y_{WV}^W(Y_{WV}^W((1+x_{2}+x_{0})^{L(0)_{n}}w, x_{0})v,
x_{2})u\nn
&&=\sum_{m = 0}^N\sum_{n = 0}^N(-1)^{m+n}{m+N \choose N}{n+N
\choose N}
\res_{x_{2}}\res_{x_{0}}\res_{x_{1}}
x_{0}^{-1}\delta\left(\frac{x_{1}-x_{2}}{x_{0}}\right)\cdot\nn
&&\quad\quad\quad\quad\cdot
x_{0}^{-N-m-1}x_{2}^{-N-n-1}
(1+x_{2}+x_{0})^{\swt w + N}(1+x_{2})^{\swt v+ N+m}
 \cdot\nn
&&\quad\quad\quad\quad\cdot
(1+x_{2}+x_{0})^{-(L(-1)+L(0))}
Y_{WV}^W((1+x_{2}+x_{0})^{L(0)_{n}}w,x_{1})Y_{V}(v, x_{2})u\nn
&&\quad -\sum_{m = 0}^N\sum_{n = 0}^N(-1)^{m+n}{m+N \choose N}{n+N
\choose N}
\res_{x_{2}}\res_{x_{0}}\res_{x_{1}}
x_{0}^{-1}\delta\left(\frac{x_{2}-x_{1}}{-x_{0}}\right)\cdot\nn
&&\quad\quad\quad\quad\cdot
x_{0}^{-N-m-1}x_{2}^{-N-n-1}
(1+x_{2}+x_{0})^{\swt w + N}(1+x_{2})^{\swt v+ N+m}
 \cdot\nn
&&\quad\quad\quad\quad\cdot
(1+x_{2}+x_{0})^{-(L(-1)+L(0))}
Y_W(v, x_{2})  Y_{WV}^W((1+x_{2}+x_{0})^{L(0)_{n}}w,x_{1})u\nn
&&=\sum_{m = 0}^N\sum_{n = 0}^N(-1)^{m+n}{m+N \choose N}{n+N
\choose N}
\res_{x_{2}}\res_{x_{0}}\res_{x_{1}}
x_{0}^{-1}\delta\left(\frac{x_{1}-x_{2}}{x_{0}}\right)\cdot\nn
&&\quad\quad\quad\quad\cdot
(x_{1}-x_{2})^{-N-m-1}x_{2}^{-N-n-1}
(1+x_{1})^{\swt w + N}(1+x_{2})^{\swt v+ N+m}
 \cdot\nn
&&\quad\quad\quad\quad\cdot 
(1+x_{1})^{-(L(-1)+L(0))}
Y_{WV}^W((1+x_{1})^{L(0)_{n}}w,x_{1})Y_{V}(v, x_{2})u\nn
&&\quad -\sum_{m = 0}^N\sum_{n = 0}^N(-1)^{m+n}{m+N \choose N}{n+N
\choose N}
\res_{x_{2}}\res_{x_{0}}\res_{x_{1}}
x_{0}^{-1}\delta\left(\frac{x_{2}-x_{1}}{-x_{0}}\right)\cdot\nn
&&\quad\quad\quad\quad\cdot
(-x_{2}+x_{1})^{-N-m-1}x_{2}^{-N-n-1}
(1+x_{1})^{\swt w + N}(1+x_{2})^{\swt v+ N+m}
 \cdot\nn
&&\quad\quad\quad\quad\cdot
(1+x_{1})^{-(L(-1)+L(0))}
Y_W(v, x_{2})  Y_{WV}^W((1+x_{1})^{L(0)_{n}}w,x_{1})u\nn
&&=\sum_{m = 0}^N\sum_{n = 0}^N
\sum_{i\in \N}(-1)^{m+n+i}{m+N \choose N}{n+N
\choose N}{-N -m - 1\choose
i}\cdot\nn
&&\quad\quad\quad\quad\cdot
\res_{x_{2}}\res_{x_{1}}
x_{1}^{-N-m-1-i}x_{2}^{-N-n-1+i}
(1+x_{2})^{m}(1+x_{1})^{-(L(-1)+L(0))}\cdot\nn
&&\quad\quad\quad\quad\cdot
Y_{WV}^W((1+x_{1})^{L(0) + N}w,x_{1})
Y_{V}((1+x_{2})^{L(0)+ N}v, x_{2})u\nn
&&\quad -\sum_{m = 0}^N\sum_{n = 0}^N
\sum_{i\in \N}(-1)^{m+n+i}{m+N \choose N}{n+N
\choose N}{-N -m - 1\choose
i}\cdot\nn
&&\quad\quad\quad\quad\cdot
\res_{x_{2}}\res_{x_{1}}
x_{1}^{i}x_{2}^{-2N-m-n-2-i}(1+x_{1})^{-(L(-1)+L(0))}\cdot\nn
&&\quad\quad\quad\quad\cdot
Y_W((1+x_{2})^{L(0)+ N+m}v, x_{2})Y_{WV}^W((1+x_{1})^{L(0)+ N}w,x_{1})u,
\end{eqnarray*}
where in the fifth step, we have used 
(\ref{a5.5}) and  Part 3 in Lemma \ref{l1} to obtain 
$(L(-1)+L(0))\tilde{w}=\omega *_{N} \tilde{w}-\tilde{w} *_{N} \omega$ which 
is in $O_{N}(W)$ for $\tilde{w}\in O_{N}(W)$ by Lemma \ref{l2},
in the sixth step, we have changed the variable $y=\frac{x_{0}}{1+x_{2}}$
and in the seventh step, we have used the Jacobi identity (\ref{jacobi-2}).
By Part 2 in Lemma  \ref{l1}, we know that every term in the second sum
in the right-hand
side is of the form $(L(-1)+L(0))^{k}\tilde{w}$ for some $k\in \N$ and 
$\tilde{w}\in O_N(W)$ and is thus in $O_N(W)$.
Also those terms with $i> N-m$ 
in the first sum in the right-hand side 
lie in $O_N(W)$. The sum of those terms with $i\le N-m$ in 
the first sum in the right-hand side
equals
\begin{eqnarray*}
\lefteqn{\sum_{m = 0}^N\sum_{n = 0}^N\sum_{i=0}^{N-m}(-1)^{m+n+i}{m+N \choose N}{n+N
\choose N}{-N -m - 1\choose
i}\cdot}\nn
&&\quad\quad\quad\quad\cdot
\res_{x_{2}}\res_{x_{1}}
x_{1}^{-N-m-1-i}x_{2}^{-N-n-1+i}
(1+x_{2})^{m}(1+x_{1})^{-(L(-1)+L(0))}\cdot\nn
&&\quad\quad\quad\quad\cdot
Y_{WV}^W((1+x_{1})^{L(0) + N}w,x_{1})
Y_{V}((1+x_{2})^{L(0)+ N}v, x_{2})u\nn
&&=w\ast_N(v\ast_N u)\nn
&&\quad +
\sum_{m = 0}^N\sum_{n = 0}^N(-1)^{m+n}{m+N \choose N}{n+N
\choose N}\res_{x_{2}}\res_{x_{1}}x_{1}^{-N-m-1}x_{2}^{-N-n-1}\cdot\nn
&&\quad\quad\quad\quad\cdot
\left(\sum_{i= 0}^{N-m}\sum_{j\in \N}{-N -m -
1\choose i}{m \choose j}(-1)^i\frac{x_{2}^{i+j}}{x_{1}^{i}} - 1\right)
\cdot\nn
&&\quad\quad\quad\quad\cdot
(1+x_{1})^{-(L(-1)+L(0))}Y_{WV}^W((1+x_{1})^{L(0)+ N}w,x_{1})
Y_{V}((1+x_{2})^{L(0)+ N}v, x_{2})u\nn
&&=w\ast_N(v\ast_N u)\nn
&&\quad +
\sum_{n = 0}^N(-1)^{n}{n+N
\choose N}\res_{x_{2}}\res_{x_{1}}x_{1}^{-N-1}x_{2}^{-N-n-1}\cdot\nn
&&\quad\quad\quad\quad\cdot
\left(\sum_{m = 0}^N(-1)^{m}{m+N \choose N}\left(\sum_{i= 0}^{N-m}\sum_{j\in \N}{-N -m -
1\choose i}{m \choose j}(-1)^i\frac{x_{2}^{i+j}}{x_{1}^{m+i}} - 
\frac{1}{x_{1}^{m}}\right)\right)
\cdot\nn
&&\quad\quad\quad\quad\cdot
(1+x_{1})^{-(L(-1)+L(0))}Y_{WV}^W((1+x_{1})^{L(0) + N}w,x_{1})
Y_{V}((1+x_{2})^{L(0)+ N}v, x_{2})u.
\end{eqnarray*}
By Proposition 5.3 in \cite{DLM},
$$\sum_{m = 0}^N(-1)^{m}{m+N \choose N}
\left(\sum_{i= 0}^{N-m}\sum_{j\in \N}{-N -m -
1\choose i}{m \choose j}(-1)^i\frac{x_{2}^{i+j}}{x_{1}^{m+i}} 
- \frac{1}{x_{1}^{m}}\right)=0.$$
Thus the calculations above give (\ref{a4}).

Finally we also have to show that the left action and the right action of
$A_N(V)$ on $A_N(W)$ commute, that is,
$$(u\ast_N w)\ast_N v \equiv u\ast_N (w \ast_N
v)\ {\rm mod}\ O_{N}(W).$$
The proof is  similar to the proof of 
(\ref{a4}) and is omitted. 
\epfv

\section{Logarithmic intertwining operators}

We begin our study of the connection between
logarithmic intertwining operators and 
associative algebras $A_{N}(V)$ in this section. Logarithmic intertwining 
operators were introduced first in the representation theory of vertex 
operator algebras by Milas \cite{M}. Here we recall the general 
definition of logarithmic intertwining operator from \cite{HLZ2}.

\begin{defn}
{\rm
Let $W_1$, $W_2$ and $W_3$ 
be generalized 
$V$-modules. A {\it logarithmic intertwining
operator of type ${W_3\choose W_1\,W_2}$} is a linear map
\begin{eqnarray*}
\mathcal{Y}: W_1\otimes W_2&\to& W_3[\log x]\{x\},
\\
w_{(1)}\otimes w_{(2)}&\mapsto &\mathcal{Y}(w_{(1)},x)w_{(2)}=\sum_{n\in
{\mathbb C}}\sum_{k\in {\mathbb N}}\mathcal{Y}_{n;\,k}(w_{(1)})
w_{(2)}x^{-n-1}(\log x)^k\in W_3[\log x]\{x\}
\end{eqnarray*}
satisfying the
following conditions:

\begin{enumerate}

\item The {\it lower truncation
condition}: For any $w_{(1)}\in W_1$, $w_{(2)}\in W_2$ and $n\in
\mathbb{C}$,
\begin{equation}\label{log:ltc}
\mathcal{Y}_{n+m;\,k}(w_{(1)})w_{(2)}=0\;\;\mbox{ for }\;m\in {\mathbb
N} \;\mbox{ sufficiently large,\, independent of}\;k.
\end{equation}

\item The {\it Jacobi identity}:
\begin{eqnarray}\label{log:jacobi}
\lefteqn{\dps x^{-1}_0\delta \bigg( \frac{x_1-x_2}{x_0}\bigg)
Y_{W_{3}}(v,x_1)\mathcal{Y}(w_{(1)},x_2)w_{(2)}}\nn
&&\hspace{2em}- x^{-1}_0\delta \bigg( \frac{x_2-x_1}{-x_0}\bigg)
\mathcal{Y}(w_{(1)},x_2)Y_{W_{2}}(v,x_1)w_{(2)}\nn
&&{\dps = x^{-1}_2\delta \bigg( \frac{x_1-x_0}{x_2}\bigg) \mathcal{
Y}(Y_{W_{1}}(v,x_0)w_{(1)},x_2) w_{(2)}}
\end{eqnarray}
for $v\in V$, $w_{(1)}\in W_1$ and $w_{(2)}\in W_2$.

\item The {\em $L(-1)$-derivative property:} for any
$w_{(1)}\in W_1$,
\begin{equation}\label{log:L(-1)dev}
\mathcal{Y}(L(-1)w_{(1)},x)=\frac d{dx}\mathcal{Y}(w_{(1)},x).
\end{equation}
\end{enumerate}}
\end{defn}

Using Proposition \ref{p2} and Corollary \ref{p2-cor}, 
we have the following result on
logarithmic intertwining operators:

\begin{prop}\label{p3}
Let $W_{1}$, $W_{2}$ and $W_{3}$ be lower-bounded generalized $V$-modules 
and $\Y$ a logarithmic  intertwining operator of type ${W_{3}\choose W_{1}W_{2}}$. 
Let $N_{2}$ and $N'_{3}$  be positive integers such that 
$W_{2}$ and $W'_{3}$ are generated by 
$\Omega^{0}_{N_{2}}(W_{2})$ and 
$\Omega^{0}_{N_{3}'}(W'_{3})$. 
For $\tilde{w}_{(1)} \in W_1$, $\tilde{w}_{(2)}
\in {W_2}$, and $\tilde{w}_{(3)}' \in W_3'$, the series $\langle
\tilde{w}_{(3)}', \mathcal {Y}(\tilde{w}_{(1)}, x)\tilde{w}_{(2)}\rangle$ can be
expressed as a linear combination of series of the form $\langle
w_{(3)}', \mathcal {Y}(w_{(1)}, x)w_{(2)}\rangle$ for $w_{(1)} \in W_1$, $w_{(2)} \in
\Omega^{0}_{N_{2}}(W_{2})$ and 
$w_{(3)}' \in \Omega^{0}_{N_{3}'}(W'_{3})$ with Laurent
polynomials of $x$  as
coefficients.
\end{prop}
\pf
By the commutator formula for vertex operators and logarithmic 
intertwining operators,
we have
$$(Y_{W_{3}})_{m}(u)\Y(w_{(1)}, x)-\Y(w_{(1)}, x)(Y_{W_{2}})_{m}(u)
=\res_{x_{0}}(x_{2}+x_{0})^{m}\Y(Y_{W_{1}}(u, x_{0})w_{(1)}, x_{2})$$
for $u\in V$ and $w_{(1)}\in W_{1}$. 
The conclusion follows from this commutator formula,  
Proposition \ref{p2}, Corollary \ref{p2-cor} and induction on the weights
of  $\tilde{w}_{(2)}$ and $\tilde{w}_{(3)}'$.
\epfv

Let
$$\mathcal{Y}^{0}(w_{(1)}, x)=\sum_{n \in \mathbb{C}} \mathcal
{Y}_{n, 0}(w_{(1)})x^{-n - 1}.$$
Then for $v\in V$, $w_{(1)}\in W_1$ and
$w_{(2)}\in W_2$, the Jacobi identity for $\mathcal {Y}^{0}(w_{(1)}, x)$ holds, that is,
\begin{eqnarray}\label{log:jacobi1}
\lefteqn{\dps x^{-1}_0\delta \bigg( \frac{x_1-x_2}{x_0}\bigg)
Y_{W_{3}}(v,x_1)\mathcal{Y}^{0}(w_{(1)},x_2)w_{(2)}}\nn
&&\hspace{2em}- x^{-1}_0\delta \bigg( \frac{x_2-x_1}{ -x_0}\bigg)
\mathcal{Y}^{0}(w_{(1)},x_2)Y_{W_{2}}(v,x_1)w_{(2)}\nn
&&{\dps = x^{-1}_2\delta \bigg( \frac{x_1-x_0}{x_2}\bigg) 
\mathcal{Y}^{0}(Y_{W_{1}}(v,x_0)w_{(1)},x_2) w_{(2)}}.
\end{eqnarray}

\begin{prop}[\cite{M}]\label{log:logwt}
Let $W_1$, $W_2$, $W_3$ be lower-bounded generalized $V$-modules
and let $\mathcal{Y}$ be a
logarithmic intertwining operator of type ${W_3\choose W_1\,W_2}$.
Let $w_{(1)}\in W_{1}$, $w_{(2)}\in W_{2}$, $h_1, h_2\in \C$,
and $k_1, k_2\in \Z_{+}$ such that
$(L(0)-h_1)^{k_1}w_{(1)}=0$ and $(L(0)-h_2)^{k_2}w_{(2)}=0$.  

\begin{enumerate}

\item For $w'_{(3)}\in W_3'$, $h_3\in \mathbb{C}$ and
$k_3\in \mathbb{Z}_+$ such that $(L'(0)-h_3)^{k_3}w'_{(3)}=0$,
\begin{eqnarray*}\label{log:k}
\lefteqn{\langle w'_{(3)}, \mathcal{Y}(w_{(1)}, x)w_{(2)}\rangle}\nn
&&\in \mathbb{C}x^{h_3-h_1-h_2}+\mathbb{C}x^{h_3-h_1-h_2}\log
x + \cdots\oplus \mathbb{C}x^{h_3-h_1-h_2}(\log
x)^{k_1+k_2+k_3-3}.\nn
\end{eqnarray*}

\item Suppose that there exist $h_3 \in {\mathbb C}$ and $k_3 \in \mathbb{Z}_+$
such that for any
homogeneous element $w'_{(3)}\in W_3'$, 
$(L'(0)- h_{3})^{k_3}w'_{(3)}=0$. Then
\begin{eqnarray*}
\mathcal{Y}(w_{(1)}, x)w_{(2)} \in
x^{h_3-h_1-h_2}W_3[[x, x^{-1}]]+ x^{h_3-h_1-h_2}W_3[[x,
x^{-1}]]\log
x\nn +  \cdots + x^{h_3-h_1-h_2}W_3[[x, x^{-1}]](\log x)^{k_1+k_2+k_3-3},\nn
\end{eqnarray*}

\end{enumerate}

\end{prop}

Let $W$ be a generalized $V$-module. 
Recall the operator $x^{\pm L(0)}$ in \cite{HLZ2} defined by
\begin{equation}\label{log:x^L(0)}
x^{\pm L(0)}w =x^{\pm n}\sum_{i\in {\mathbb
N}}\frac{(L(0)-n)^iw}{i!}(\pm\log x)^i\in x^{\pm n}W_{[n]}[\log x]
\end{equation}
for $w\in W_{[n]}$. Also recall the $L(0)$-conjugation property
for logarithmic intertwining operators in \cite{HLZ2}: For any logarithmic 
intertwining operator $\mathcal{Y}$ of type ${W_{3}\choose W_{1}W_{2}}$ and any
$w_{(1)}\in W_{1}$,
\begin{equation}\label{l-0-conj-formal}
y^{L(0)}\mathcal{Y}(w_{(1)}, x)y^{-L(0)}=\mathcal{Y}(y^{L(0)}w_{(1)}, xy).
\end{equation}


%
%

Let $W_1, W_2$ and $W_3$ be lower-bounded 
generalized $V$-modules and $\mathcal
{Y}$ an logarithmic intertwining operator of type ${W_3\choose
W_1\,W_2}$. Then for $N\in \N$, $A_N(W_1)
\otimes_{A_N(V)} \Omega^{0}_N(W_2)$ and $\Omega^{0}_N(W_3)$ are both left
$A_N(V)$-modules. Note that $\Omega^{0}_N(W_3)$ as a subspace of $W_{3}$ is also 
graded and for any $n\in \C$, the image of $\Omega^{0}_N(W_3)$ under
the projection $\pi_{n}: W_{3}\to (W_{3})_{[n]}$ is denoted
$(\Omega^{0}_N(W_3))_{[n]}$.

First, we consider the case that there exists $h_{3}\in \C$ 
such that $W_{3}=\coprod_{n\in h_{3}+\N}(W_{3})_{[n]}$ and 
$(W_{3})_{[h_{3}]}\ne 0$.
As above, let $L(0)_{s}$ be the semisimple part of the operator $L(0)$
on any module for the Virasoro algebra. Let 
\[
\rho(\mathcal {Y}): W_1\otimes \Omega^{0}_N(W_2)
\rightarrow W_{3}
\]
be defined by
\begin{eqnarray*}
\rho(\mathcal {Y})(w_{(1)}\otimes w_{(2)}) 
&=&\sum_{n=0}^{N}{\rm
Res}_x x^{-h_3 -n-
1}\mathcal{Y}^{0}(x^{L(0)_{s}}w_{(1)},
x)x^{L(0)_{s}}w_{(2)}\nn
& =& \sum_{n=0}^{N}\Y_{\swt w_{(1)} + \swt w_{(2)} -
h_3 -n - 1, 0}(w_{(1)})w_{(2)},\nn
\end{eqnarray*}
for  homogeneous $w_{(1)} \in W_1$, $w_{(2)} \in \Omega^{0}_N(W_2)$. 
In the general case that 
$W_{3}=\coprod_{\mu\in \C/\Z}W_{3}^{\mu}$,
let $\pi^{\mu}: W_{3}\to W_{3}^{\mu}$ be the projection from
$W_{3}$ to $W_{3}^{\mu}$ for $\mu\in \C/\Z$. Then 
$\pi^{\mu}\circ \mathcal{Y}$ is a logarithmic intertwining operator
of type ${W_3^{\mu}\choose
W_1\,W_2}$ for each $\mu\in \C/\Z$. 
We define 
\[
\rho(\mathcal {Y}): W_1\otimes \Omega^{0}_N(W_2)
\rightarrow W_{3}
\]
by 
$$\rho(\mathcal {Y})(w_{(1)}\otimes w_{(2)}) = \sum_{\mu\in \C/\Z}
\rho(\pi^{\mu}\circ \mathcal {Y})(w_{(1)}\otimes w_{(2)})$$
for $w_{(1)}\in W_{1}$ and $w_{(2)}\in \Omega^{0}_N(W_2)$.
We have:

\begin{lemma}
The image of $W_1\otimes \Omega^{0}_N(W_2)$ under $\rho(\Y)$ 
is in $\Omega^{0}_N(W_3)$. In particular, $\rho(\Y)$ is a linear map from 
$W_1\otimes \Omega^{0}_N(W_2)$
to $\Omega^{0}_N(W_3)$.
\end{lemma}
\pf
This follows from the definition by calculating the weights. 
\epfv

The following lemma shows that $\rho(\mathcal
{Y})$ is in fact a linear map from  $A_{N}(W_1)\otimes \Omega^{0}_N(W_2)$
to $\Omega^{0}_N(W_3)$:

\begin{lemma}\label{rho-1}
For $w_{(1)} \in O_N(W_1)$ and $w_{(2)} \in \Omega^{0}_N(W_2)$, 
$\rho(\mathcal {Y})(w_{(1)}\otimes
w_{(2)}) = 0$.
\end{lemma} 
\pf 
We prove the lemma in the case that
$W_{3}=\coprod_{n\in h_{3}+\N}(W_{3})_{[n]}$ for some $h_{3}\in \C$ and
$(W_{3})_{[h_{3}]}\ne 0$,
\begin{eqnarray}\label{w-1}
w_{(1)} &=&u\circ_{N} w\nn
&=&
{\rm Res}_x x^{-2N-2}Y_{W_{1}}((1+x)^{L(0)+N}u, x)w\nn
&=&
{\rm Res}_x x^{-2N-2}(1+x)^{\swt u+N}Y_{W_{1}}(u, x)w
\end{eqnarray}
for some homogeneous $u \in V$ and $w \in W_1$,
and $w_{(2)} \in (\Omega^{0}_N(W_2))_{[h_2]}$.
The general case follows easily. 

In this case, by the definition of $\rho(\Y)$, (\ref{w-1}),
the $L(0)_{s}$-conjugation property above
and the Jacobi identity for $\Y^{0}$,
we have
\begin{eqnarray}\label{rho-y-1}
\lefteqn{\rho(\mathcal {Y})((w_{(1)} + O_N(W_1))\otimes w_{(2)})}\nn
&& = \sum_{n=0}^{N}  \mbox{\rm Res}_{x_2} \res_{x_{0}}
x_{0}^{-2N -
2}x_{2}^{h_{2}-h_{3}-n-1}(1+x_{0})^{N}\cdot\nn
&&\quad\quad\quad\quad \quad\quad\quad\quad \cdot
\mathcal{Y}^{0}(x_{2}^{L(0)_{s}}Y_{W_{1}}((1+x_{0})^{L(0)}u, x_{0})w_{(1)}, x_2)w_{(2)}
\nn
&& = \sum_{n=0}^{N}  \mbox{\rm Res}_{x_2} \res_{x_{0}}
x_{0}^{-2N -
2}x_{2}^{h_{2}-h_{3}-n-1}(1+x_{0})^{N}\cdot\nn
&&\quad\quad\quad\quad \quad\quad\quad\quad  \cdot
\mathcal{Y}^{0}(Y_{W_{1}}((x_{2}+x_{0}x_{2})^{L(0)_{s}}u, x_{0}x_{2})
x_{2}^{L(0)_{s}}w_{(1)}, x_2)w_{(2)}\nn
&& = \sum_{n=0}^{N}  \mbox{\rm Res}_{x_2} \res_{x_{0}}\res_{x_{1}}
x_{1}^{-1}\delta\left(\frac{x_{2}+x_{0}x_{2}}{x_{1}}\right)x_{0}^{-2N -
2}x_{2}^{h_{2}-h_{3}-n-N-1}(x_{2}+x_{0}x_{2})^{\swt u+N}\cdot\nn
&&\quad\quad\quad\quad \quad\quad\quad\quad  \cdot
\mathcal{Y}^{0}(Y_{W_{1}}(u, x_{0}x_{2})
x_{2}^{L(0)_{s}}w_{(1)}, x_2)w_{(2)}\nn
&& = \sum_{n=0}^{N}  \mbox{\rm Res}_{x_2} \res_{x_{0}}\res_{x_{1}}
x_{1}^{-1}\delta\left(\frac{x_{2}+x_{0}x_{2}}{x_{1}}\right)x_{0}^{-2N -
2}x_{2}^{h_{2}-h_{3}-n-N-1}x_{1}^{\swt u+N}\cdot\nn
&&\quad\quad\quad\quad \quad\quad\quad\quad  \cdot
\mathcal{Y}^{0}(Y_{W_{1}}(u, x_{0}x_{2})
x_{2}^{L(0)_{s}}w_{(1)}, x_2)w_{(2)}\nn
&&=\sum_{n=0}^{N}  \mbox{\rm Res}_{x_2} \res_{x_{0}}\res_{x_{1}}
(x_{0}x_{2})^{-1}\delta\left(\frac{x_{1}-x_{2}}{x_{0}x_{2}}\right)
x_{0}^{-2N -
2}x_{2}^{h_{2}-h_{3}-n-N-1}\cdot\nn
&&\quad\quad\quad\quad \quad\quad\quad\quad \cdot
x_{1}^{\swt u+N}Y_{W_{3}}(u, x_{1})
\mathcal{Y}^{0}(x_{2}^{L(0)_{s}}w_{(1)}, x_2)w_{(2)}
\nn
&&\quad -\sum_{n=0}^{N}  \mbox{\rm Res}_{x_2} \res_{x_{0}}\res_{x_{1}}
(x_{0}x_{2})^{-1}\delta\left(\frac{x_{2}-x_{1}}{-x_{0}x_{2}}\right)
x_{0}^{-2N -
2}x_{2}^{h_{2}-h_{3}-n-N-1}\cdot\nn
&&\quad\quad\quad\quad \quad\quad\quad\quad \cdot x_{1}^{\swt u+N}
\mathcal{Y}^{0}(x_{2}^{L(0)_{s}}w_{(1)}, x_2)
Y_{W_{2}}(u, x_{1})w_{(2)}\nn
&&=\sum_{n=0}^{N} 
\mbox{\rm Res}_{x_2} \res_{x_{1}}(x_{1}-x_{2})^{-2N -
2}x_{2}^{h_{2}-h_{3}-n+N}\cdot\nn
&&\quad\quad\quad\quad \quad\quad\quad\quad \cdot
x_{1}^{\swt u+N}Y_{W_{3}}(u, x_{1})
\mathcal{Y}^{0}(x_{2}^{L(0)_{s}}w_{(1)}, x_2)w_{(2)}\nn
&&\quad -\sum_{n=0}^{N} 
\mbox{\rm Res}_{x_2} \res_{x_{1}}(-x_{2}+x_{1})^{-2N -
2}x_{2}^{h_{2}-h_{3}-n+N}\cdot\nn
&&\quad\quad\quad\quad \quad\quad\quad\quad \cdot 
x_{1}^{\swt u+N}\mathcal{Y}^{0}(x_{2}^{L(0)_{s}}w_{(1)}, x_2)
Y_{W_{2}}(u, x_{1})w_{(2)}
\nn
&&=\sum_{n=0}^{N} \sum_{m = 0}^N  {-2N-2\choose l}
x_{2}^{h_{2}-h_{3}-n+N+l}\cdot\nn
&&\quad\quad\quad\quad \quad\quad\quad\quad \cdot
\mbox{\rm Res}_{x_2} (Y_{W_{3}})_{\swt u-N-2-l}(u)
\mathcal{Y}^{0}(x_{2}^{L(0)_{s}}w_{(1)}, x_2)w_{(2)}\nn
&&\quad -\sum_{n=0}^{N}  \sum_{l\in \N}(-1)^{m-2N-2-l}  {-2N-2\choose l}
x_{2}^{h_{2}-h_{3}-n-N-2-l}\cdot\nn
&&\quad\quad\quad\quad \quad\quad\quad\quad \cdot 
\mbox{\rm Res}_{x_2} \mathcal{Y}^{0}(x_{2}^{L(0)_{s}}w_{(1)}, x_2)
(Y_{W_{2}})_{\swt u +N+l}(u)w_{(2)}.
\end{eqnarray}
Since the weight of $(Y_{W_{3}})_{\swt u-N-2-l}(u)$ 
is $N+1+l$, the real parts of the weights of 
the homogeneous components of 
the coefficients of 
$$(Y_{W_{3}})_{\swt u-N-2-l}(u)
\mathcal{Y}^{0}(x_{2}^{L(0)_{s}}w_{(1)}, x_2)w_{(2)}$$
are larger than or equal to $N+1+l+r_{3}$. Thus the first term in the right-hand side
of (\ref{rho-y-1}) is equal to $0$. Since $w_{(2)}\in \Omega^{0}_N(W_2)$, 
$(Y_{W_{2}})_{\swt u +N+l}(u)w_{(2)}=0$. From (\ref{rho-y-1}), 
we see that the second term in the right-hand side of 
(\ref{rho-y-1}) is also equal to $0$. Thus 
the left-hand side of (\ref{rho-y-1})
is $0$, 
proving the lemma.
\epfv

As we mentioned before, by this lemma, $\rho(\Y)$ is in fact 
a linear map from $A_{N}(W_1)\otimes \Omega^{0}_N(W_2)$
to $\Omega^{0}_N(W_3)$. 
We now have:

\begin{prop}
The map $\rho(\mathcal {Y})$ is in fact an $A_N(V)$-module
homomorphism from $A_N(W_1)\otimes \Omega^{0}_N(W_2)$ to 
$\Omega^{0}_N(W_3)$, that is, 
$$\rho(\mathcal {Y}) \in {\rm Hom}_{A_N(V)}
(A_N(W_1)\otimes \Omega^{0}_N(W_2), \Omega^{0}_N(W_3)).$$
\end{prop}
\pf 
We need to prove 
\[
\rho(\mathcal {Y})((u *_N w_{(1)} + O_N(W_1)) \otimes
w_{(2)}) = o(u) \rho(\mathcal {Y})((w_{(1)} + O_N(W_1)) \otimes w_{(2)}).
\]
for $u \in V$, $w_{(1)} \in W_1$, and $w_{(2)} \in \Omega^{0}_N(W_2)$.
We prove this in the case that 
$W_{3}=\coprod_{n\in h_{3}+\N}(W_{3})_{[n]}$ for some $h_{3}\in \C$ and
$(W_{3})_{[h_{3}]}\ne 0$. The general case follows easily. 

Let $w_{(2)} \in \Omega^{0}_N(W_2)$ be homogeneous
of weight $h_{2}$. 
Calculations similar to those in the proof of Lemma \ref{rho-1} give
\begin{eqnarray}\label{rho-y-homo-1}
\lefteqn{\rho(\mathcal {Y})((u *_N w_{(1)} + O_N(W_1))\otimes w_{(2)})}\nn
&& = \sum_{n=0}^{N} \sum_{m = 0}^N (-1)^m {m+N\choose
N}\mbox{\rm Res}_{x_2} \res_{x_{0}}
 x_{0}^{-N -m-1}x_{2}^{h_{2}-h_{3}-n-1}(1+x_{0})^{N}\cdot\nn
&&\quad\quad\quad\quad\quad\quad\quad\quad  \cdot
\mathcal{Y}^{0}(x_{2}^{L(0)_{s}}Y_{W_{1}}
((1+x_{0})^{L(0)}u, x_{0})w_{(1)}, x_2)w_{(2)}
\nn
&& = \sum_{n=0}^{N} \sum_{m = 0}^N (-1)^m {m+N\choose
N} \mbox{\rm Res}_{x_2} \res_{x_{0}}
x_{0}^{-N -m-1}x_{2}^{h_{2}-h_{3}-n-1}(1+x_{0})^{N}\cdot\nn
&&\quad\quad\quad\quad\quad\quad\quad\quad  \cdot
\mathcal{Y}^{0}(Y_{W_{1}}((x_{2}+x_{0}x_{2})^{L(0)}u, x_{0}x_{2})
x_{2}^{L(0)_{s}}w_{(1)}, x_2)w_{(2)}
\nn
&& = \sum_{n=0}^{N} \sum_{m = 0}^N (-1)^m {m+N\choose
N} \mbox{\rm Res}_{x_2} \res_{x_{0}}\res_{x_{1}}
x_{1}^{-1}\delta\left(\frac{x_{2}+x_{0}x_{2}}{x_{1}}\right)
x_{0}^{-N -m-1}\cdot\nn
&&\quad\quad\quad\quad\quad\quad\quad\quad  \cdot
x_{2}^{h_{2}-h_{3}-n-N-1}(x_{2}+x_{0}x_{2})^{\swt u+N}
\mathcal{Y}^{0}(Y_{W_{1}}(u, x_{0}x_{2})
x_{2}^{L(0)_{s}}w_{(1)}, x_2)w_{(2)}
\nn
&& = \sum_{n=0}^{N} \sum_{m = 0}^N (-1)^m {m+N\choose
N} \mbox{\rm Res}_{x_2} \res_{x_{0}}\res_{x_{1}}
x_{1}^{-1}\delta\left(\frac{x_{2}+x_{0}x_{2}}{x_{1}}\right)
x_{0}^{-N -m-1}\cdot\nn
&&\quad\quad\quad\quad\quad\quad\quad\quad  \cdot
x_{2}^{h_{2}-h_{3}-n-N-1}x_{1}^{\swt u+N}
\mathcal{Y}^{0}(Y_{W_{1}}(u, x_{0}x_{2})
x_{2}^{L(0)_{s}}w_{(1)}, x_2)w_{(2)}
\nn
&& =\sum_{n=0}^{N} \sum_{m = 0}^N (-1)^m {m+N\choose
N} \mbox{\rm Res}_{x_2} \res_{x_{0}}\res_{x_{1}}
(x_{0}x_{2})^{-1}\delta\left(\frac{x_{1}-x_{2}}{x_{0}x_{2}}\right)
\cdot\nn
&&\quad\quad\quad\quad \quad\quad\quad\quad \cdot
x_{2}^{h_{2}-h_{3}-n-N-1}x_{0}^{-N -m-
1}x_{1}^{\swt u+N}
Y_{W_{3}}(u, x_{1})
\mathcal{Y}^{0}(x_{2}^{L(0)_{s}}w_{(1)}, x_2)w_{(2)}
\nn
&&\quad -\sum_{n=0}^{N} \sum_{m = 0}^N (-1)^m {m+N\choose
N} \mbox{\rm Res}_{x_2} \res_{x_{0}}\res_{x_{1}}
(x_{0}x_{2})^{-1}\delta\left(\frac{x_{2}-x_{1}}{-x_{0}x_{2}}\right)
\cdot\nn
&&\quad\quad\quad\quad\quad\quad\quad\quad  \cdot 
x_{2}^{h_{2}-h_{3}-n-N-1}x_{0}^{-N-m-
1}x_{1}^{\swt u+N}
\mathcal{Y}^{0}(x_{2}^{L(0)_{s}}w_{(1)}, x_2)
Y_{W_{2}}(u, x_{1})w_{(2)}\nn
&& =\sum_{n=0}^{N} \sum_{m = 0}^N (-1)^m {m+N\choose
N} \mbox{\rm Res}_{x_2} \res_{x_{1}}
x_{2}^{h_{2}-h_{3}-n+m-1}(x_{1}-x_{2})^{-N -m-
1}x_{1}^{\swt u+N}\cdot\nn
&&\quad\quad\quad\quad \quad\quad\quad\quad \cdot
Y_{W_{3}}(u, x_{1})
\mathcal{Y}^{0}(x_{2}^{L(0)_{s}}w_{(1)}, x_2)w_{(2)}
\nn
&&\quad -\sum_{n=0}^{N} \sum_{m = 0}^N (-1)^m {m+N\choose
N} \mbox{\rm Res}_{x_2} \res_{x_{1}}
x_{2}^{h_{2}-h_{3}-n+m-1}(-x_{2}+x_{1})^{-N-m-
1}x_{1}^{\swt u+N}\cdot\nn
&&\quad\quad\quad\quad\quad\quad\quad\quad  \cdot 
\mathcal{Y}^{0}(x_{2}^{L(0)_{s}}w_{(1)}, x_2)
Y_{W_{2}}(u, x_{1})w_{(2)}\nn
&&=\sum_{n=0}^{N} \sum_{m = 0}^N \sum_{l\in \N}(-1)^{m+l} {m+N\choose
N} {-N-m-1\choose l}
\mbox{\rm Res}_{x_2}x_{2}^{h_{2}-h_{3}-n+m+l-1}\cdot\nn
&&\quad\quad\quad\quad \quad\quad\quad\quad  \cdot
 (Y_{W_{3}})_{\swt u-m-1-l}(u)
\mathcal{Y}^{0}(x_{2}^{L(0)_{s}}w_{(1)}, x_2)w_{(2)}
\nn
&&\quad -\sum_{n=0}^{N} \sum_{m = 0}^N \sum_{l\in \N}(-1)^{m-N-m-1-l} {m+N\choose
N} {-N-m-1\choose l}
\mbox{\rm Res}_{x_2}x_{2}^{h_{2}-h_{3}-n-N-l-2}\cdot\nn
&&\quad\quad\quad\quad\quad\quad\quad\quad   \cdot 
 \mathcal{Y}^{0}(x_{2}^{L(0)_{s}}w_{(1)}, x_2)
(Y_{W_{2}})_{\swt u +N+l}(u)w_{(2)}
\nn
&&=\sum_{n=0}^{N} \sum_{m = 0}^N \sum_{l=0}^{N-m}(-1)^{m+l} {m+N\choose
N} {-N-m-1\choose l}
\mbox{\rm Res}_{x_2}x_{2}^{h_{2}-h_{3}-n+m+l-1}\cdot\nn
&&\quad\quad\quad\quad \quad\quad\quad\quad  \cdot
 (Y_{W_{3}})_{\swt u-m-1-l}(u)
\mathcal{Y}^{0}(x_{2}^{L(0)_{s}}w_{(1)}, x_2)w_{(2)}
\nn
&&=\sum_{n=0}^{N} \sum_{m = 0}^N \sum_{i=m}^{N}(-1)^{i} {m+N\choose
N} {-N-m-1\choose i-m}
\mbox{\rm Res}_{x_2}x_{2}^{h_{2}-h_{3}-n+i-1}\cdot\nn
&&\quad\quad\quad\quad\quad\quad\quad\quad  \cdot
(Y_{W_{3}})_{\swt u-i-1}(u)
\mathcal{Y}^{0}(x_{2}^{L(0)_{s}}w_{(1)}, x_2)w_{(2)}
\nn
&&=\sum_{n=0}^{N} \sum_{i = 0}^N \sum_{m=0}^{i}(-1)^{i} {m+N\choose
N} {-N-m-1\choose i-m}
\mbox{\rm Res}_{x_2} x_{2}^{h_{2}-h_{3}-n+i-1}\cdot\nn
&&\quad\quad\quad\quad\quad\quad\quad\quad  \cdot
(Y_{W_{3}})_{\swt u-i-1}(u)
\mathcal{Y}^{0}(x_{2}^{L(0)_{s}}w_{(1)}, x_2)w_{(2)}
\nn
&&=\sum_{n=0}^{N}  
\mbox{\rm Res}_{x_2} x_{2}^{h_{2}-h_{3}-n-1}(Y_{W_{3}})_{\swt u-1}(u)
\mathcal{Y}^{0}(x_{2}^{L(0)_{s}}w_{(1)}, x_2)w_{(2)}
\nn
&&\quad +\sum_{n=0}^{N} \sum_{i = 1}^N \sum_{m=0}^{i}(-1)^{i} {m+N\choose
N} {-N-m-1\choose i-m}
x_{2}^{h_{2}-h_{3}-n+i-1}\cdot\nn
&&\quad\quad\quad\quad\quad\quad\quad\quad  \cdot
\mbox{\rm Res}_{x_2} (Y_{W_{3}})_{\swt u-i-1}(u)
\mathcal{Y}^{0}(x_{2}^{L(0)_{s}}w_{(1)}, x_2)w_{(2)}.
\end{eqnarray}

Since for $i=1, \dots, N$,
\begin{eqnarray}\label{comb-ident}
\lefteqn{\sum_{m = 0}^{i}{m+N\choose
N} {-N-m-1\choose i-m}}\nn
&& = \sum_{m = 0}^{i}(-1)^{i - m}{m+N\choose
N} {N+i\choose i-m}\nn
&& = \sum_{m = 0}^{i}(-1)^{i - m}{N+i\choose
N} {i\choose m}\nn
&& = {N+i\choose
N} (1-1)^{i}\nn
&&= 0,
\end{eqnarray}
the right-hand side of (\ref{rho-y-homo-1}) is equal to 
\begin{eqnarray*} 
\lefteqn{\sum_{n=0}^{N}  
\mbox{\rm Res}_{x_2} (Y_{W_{3}})_{\swt u-1}(u)
\mathcal{Y}^{0}(x_{2}^{L(0)_{s}}w_{(1)}, x_2)w_{(2)}
x_{2}^{h_{2}-h_{3}-n-1}}\nn
&& \;\;\;\;\;\;\;\;\;\;\;\;\;\;\;\;\;
= o(u) \rho(\mathcal {Y})((w_{(1)} + O_N(W_1)) \otimes w_{(2)}).
\end{eqnarray*}
This completes the proof.
\epfv

\begin{prop}\label{over-a-n}
The map $\rho(\mathcal {Y})$ is in fact an $A_N(V)$-module
homomorphism from $A_N(W_1)\otimes_{A_{N}(V)} \Omega^{0}_N(W_2)$ to 
$\Omega^{0}_N(W_3)$, that is, 
$$\rho(\mathcal {Y}) \in {\rm Hom}_{A_N(V)}
(A_N(W_1)\otimes_{A_{N}(V)} \Omega^{0}_N(W_2), \Omega^{0}_N(W_3)).$$
\end{prop}
\pf 
We need to prove 
$$\rho(\mathcal {Y})((w_{(1)} *_N u + O_N(W_1)) \otimes
w_{(2)})= \rho(\mathcal {Y})((w_{(1)} + O_N(W_1)) \otimes o(u) w_{(2)})$$
for $u \in V$, $w_{(1)} \in W_1$, and $w_{(2)} \in \Omega^{0}_N(W_2)$.
We shall prove this equality  only in the case that 
$W_{3}=\coprod_{n\in h_{3}+\N}(W_{3})_{[n]}$ for some $h_{3}\in \C$ and
$(W_{3})_{[h_{3}]}\ne 0$. The general case follows easily.

Let $w_{(2)} \in \Omega^{0}_N(W_2)$ be homogeneous of weight $h_{2}$.
By Part 1 in Lemma \ref{l1},
\begin{eqnarray*}
\lefteqn{w_{(1)} *_N u+ O_N(W_1)}\nn
&&=\sum_{m = 0}^N {m+N\choose N}(-1)^N {\rm Res}_x 
x^{-N-m-1} Y_W((1+ x)^{L(0)+m-1}u, x)w_{(1)}+ O_N(W_1).
\end{eqnarray*}
Then 
\begin{eqnarray}\label{rho-y-an-1}
\lefteqn{\rho(\mathcal {Y})((w_{(1)} *_N u + O_N(W_1)) \otimes
w_{(2)})}\nn
&& = \sum_{n=0}^{N} \sum_{m = 0}^N (-1)^N {m+N\choose
N} \mbox{\rm Res}_{x_2} \res_{x_{0}}
x_{0}^{-N -m-1}x_{2}^{h_{2}-h_{3}-n-1}(1+x_{0})^{m-1}\cdot\nn
&&\quad\quad\quad\quad\quad\quad\quad\quad  \cdot
\mathcal{Y}^{0}(x_{2}^{L(0)_{s}}Y_{W_{1}}
((1+x_{0})^{L(0)}u, x_{0})w_{(1)}, x_2)w_{(2)}
\nn
&& = \sum_{n=0}^{N} \sum_{m = 0}^N (-1)^N {m+N\choose
N} \mbox{\rm Res}_{x_2} \res_{x_{0}}
x_{0}^{-N -m-1}x_{2}^{h_{2}-h_{3}-n-1}(1+x_{0})^{m-1}\cdot\nn
&&\quad\quad\quad\quad\quad\quad\quad\quad  \cdot
\mathcal{Y}^{0}(Y_{W_{1}}((x_{2}+x_{0}x_{2})^{L(0)}u, x_{0}x_{2})
x_{2}^{L(0)_{s}}w_{(1)}, x_2)w_{(2)}
\nn
&& = \sum_{n=0}^{N} \sum_{m = 0}^N (-1)^N {m+N\choose
N} \mbox{\rm Res}_{x_2} \res_{x_{0}}\res_{x_{1}}
x_{1}^{-1}\delta\left(\frac{x_{2}+x_{0}x_{2}}{x_{1}}\right)
x_{0}^{-N -m-1}x_{2}^{h_{2}-h_{3}-n-m}\cdot\nn
&&\quad\quad\quad\quad\quad\quad\quad\quad  \cdot
(x_{2}+x_{0}x_{2})^{\swt u+m-1}\mathcal{Y}^{0}(Y_{W_{1}}(u, x_{0}x_{2})
x_{2}^{L(0)_{s}}w_{(1)}, x_2)w_{(2)}
\nn
&& = \sum_{n=0}^{N} \sum_{m = 0}^N (-1)^N {m+N\choose
N} \mbox{\rm Res}_{x_2} \res_{x_{0}}\res_{x_{1}}
x_{1}^{-1}\delta\left(\frac{x_{2}+x_{0}x_{2}}{x_{1}}\right)
x_{0}^{-N -m-1}x_{2}^{h_{2}-h_{3}-n-m}\cdot\nn
&&\quad\quad\quad\quad\quad\quad\quad\quad  \cdot
x_{1}^{\swt u+m-1}\mathcal{Y}^{0}(Y_{W_{1}}(u, x_{0}x_{2})
x_{2}^{L(0)_{s}}w_{(1)}, x_2)w_{(2)}
\nn
&&=\sum_{n=0}^{N} \sum_{m = 0}^N (-1)^N {m+N\choose
N} \mbox{\rm Res}_{x_2} \res_{x_{0}}\res_{x_{1}}
(x_{0}x_{2})^{-1}\delta\left(\frac{x_{1}-x_{2}}{x_{0}x_{2}}\right)
\cdot\nn
&&\quad\quad\quad\quad \quad\quad\quad\quad \cdot
x_{2}^{\swt w_{(1)}+h_{2}-h_{3}-n-m}x_{0}^{-N -m-
1}x_{1}^{\swt u+m-1}
Y_{W_{3}}(u, x_{1})
\mathcal{Y}^{0}(w_{(1)}, x_2)w_{(2)}
\nn
&&\quad -\sum_{n=0}^{N} \sum_{m = 0}^N (-1)^N {m+N\choose
N} \mbox{\rm Res}_{x_2} \res_{x_{0}}\res_{x_{1}}
(x_{0}x_{2})^{-1}\delta\left(\frac{x_{2}-x_{1}}{-x_{0}x_{2}}\right)
\cdot\nn
&&\quad\quad\quad\quad\quad\quad\quad\quad  \cdot 
x_{2}^{h_{2}-h_{3}-n-m}x_{0}^{-N-m-
1}x_{1}^{\swt u+m-1}
\mathcal{Y}^{0}(x_{2}^{L(0)_{s}}w_{(1)}, x_2)
Y_{W_{2}}(u, x_{1})w_{(2)}\nn
&& =\sum_{n=0}^{N} \sum_{m = 0}^N (-1)^N {m+N\choose
N} \mbox{\rm Res}_{x_2} \res_{x_{1}}
x_{2}^{\swt w_{(1)}+h_{2}-h_{3}-n+N}(x_{1}-x_{2})^{-N -m-
1}\cdot\nn
&&\quad\quad\quad\quad \quad\quad\quad\quad \cdot
x_{1}^{\swt u+m-1}Y_{W_{3}}(u, x_{1})
\mathcal{Y}^{0}(w_{(1)}, x_2)w_{(2)}
\nn
&&\quad -\sum_{n=0}^{N} \sum_{m = 0}^N (-1)^N {m+N\choose
N} \mbox{\rm Res}_{x_2} \res_{x_{1}}
x_{2}^{h_{2}-h_{3}-n+N}(-x_{2}+x_{1})^{-N-m-
1}\cdot\nn
&&\quad\quad\quad\quad\quad\quad\quad\quad  \cdot 
x_{1}^{\swt u+m-1}\mathcal{Y}^{0}(x_{2}^{L(0)_{s}}w_{(1)}, x_2)
Y_{W_{2}}(u, x_{1})w_{(2)}\nn
&&=\sum_{n=0}^{N} \sum_{m = 0}^N \sum_{l\in \N}(-1)^{N+l} {m+N\choose
N} {-N-m-1\choose l}\cdot\nn
&&\quad\quad\quad\quad \quad\quad\quad\quad  \cdot
\mbox{\rm Res}_{x_1} x_{1}^{\swt u-N-2-l}Y_{W_{3}}(u, x_{1})
\mathcal{Y}_{\swt w_{(1)}+h_{2}-h_{3}-n+N+l, 0}(w_{(1)})w_{(2)}
\nn
&&\quad -\sum_{n=0}^{N} \sum_{m = 0}^N \sum_{l\in \N}(-1)^{-m-1-l} {m+N\choose
N} {-N-m-1\choose l}
\cdot\nn
&&\quad\quad\quad\quad\quad\quad\quad\quad   \cdot 
\mbox{\rm Res}_{x_2} x_{2}^{h_{2}-h_{3}-n-m-l-1}
\mathcal{Y}^{0}(x_{2}^{L(0)_{s}}w_{(1)}, x_{2})
(Y_{W_{2}})_{\swt u+m-1+l}(u)w_{(2)}
\nn
&&=\sum_{n=0}^{N} \sum_{m = 0}^N \sum_{l=0}^{N-m}(-1)^{m+l} {m+N\choose
N} {-N-m-1\choose l}
\cdot\nn
&&\quad\quad\quad\quad \quad\quad\quad\quad  \cdot
\mbox{\rm Res}_{x_2}x_{2}^{h_{2}-h_{3}-n-m-l-1} \mathcal{Y}^{0}(x_{2}^{L(0)_{s}}w_{(1)}, x_{2})
(Y_{W_{2}})_{\swt u+m-1+l}(u)w_{(2)}
\nn
&&=\sum_{n=0}^{N} \sum_{m = 0}^N \sum_{i=m}^{N}(-1)^{i} {m+N\choose
N} {-N-m-1\choose i-m}
\cdot\nn
&&\quad\quad\quad\quad\quad\quad\quad\quad  \cdot
\mbox{\rm Res}_{x_2} x_{2}^{h_{2}-h_{3}-n-i-1}
\mathcal{Y}^{0}(x_{2}^{L(0)_{s}}w_{(1)}, x_{2})
(Y_{W_{2}})_{\swt u+i-1}(u)w_{(2)}
\nn
&&=\sum_{n=0}^{N} \sum_{i = 0}^N \sum_{m=0}^{i}(-1)^{i} {m+N\choose
N} {-N-m-1\choose i-m}
\cdot\nn
&&\quad\quad\quad\quad\quad\quad\quad\quad  \cdot
\mbox{\rm Res}_{x_2} x_{2}^{h_{2}-h_{3}-n-i-1}
\mathcal{Y}^{0}(x_{2}^{L(0)_{s}}w_{(1)}, x_{2})
(Y_{W_{2}})_{\swt u+i-1}(u)w_{(2)}
\nn
&&=\sum_{n=0}^{N}  
\mbox{\rm Res}_{x_2} x_{2}^{h_{2}-h_{3}-n-1}\mathcal{Y}^{0}(x_{2}^{L(0)_{s}}w_{(1)}, x_{2})
(Y_{W_{2}})_{\swt u-1}(u)w_{(2)}
\nn
&&\quad +\sum_{n=0}^{N} \sum_{i = 1}^N \sum_{m=0}^{i}(-1)^{i} {m+N\choose
N} {-N-m-1\choose i-m}
\cdot\nn
&&\quad\quad\quad\quad\quad\quad\quad\quad  \cdot
\mbox{\rm Res}_{x_2} x_{2}^{h_{2}-h_{3}-n-i-1}\mathcal{Y}^{0}(x_{2}^{L(0)_{s}}w_{(1)}, x_{2})
(Y_{W_{2}})_{\swt u+i-1}(u)w_{(2)}\nn
&&=\rho(\mathcal {Y})((w_{(1)} + O_N(W_1)) \otimes o(u) w_{(2)}),
\end{eqnarray}
where in the last step, we used (\ref{comb-ident}).
This completes the proof.
\epfv

Let  $\mathcal{Y}$ be a logarithmic intertwining operator of type 
${W_{3}\choose W_{1} W_{2}}$. For homogeneous 
$w_{(1)}\in W_{1}$ and $w_{(2)}\in W_{2}$
and $w_{(3)}'\in W_{3}'$, let $k_{1}, k_{2}, k_{3}\in \Z_{+}$ such that 
$(L(0)-\wt w_{(1)})^{k_{1}}w_{(1)}=0$, 
$(L(0)-\wt w_{(2)})^{k_{2}}w_{(2)}=0$ and 
$(L'(0)-\wt w'_{(3)})^{k_{3}}w'_{(3)}=0$.
Then by Proposition \ref{log:logwt},  
\begin{eqnarray*}
\lefteqn{\langle w'_{(3)}, \mathcal{Y}(w_{(1)}, x)w_{(2)}\rangle}\nn
&&\in 
\C x^{\swt w_{(3)}'-\swt w_{(1)}-\swt w_{(2)}}+
\C x^{\swt w_{(3)}'-\swt w_{(1)}-\swt w_{(2)}}\log x\nn
&&\quad +\cdots +\C x^{\swt w_{(3)}'-\swt w_{(1)}-\swt w_{(2)}}
(\log x)^{k_{1}+k_{2}+k_{3}-3}.
\end{eqnarray*}
For $z\in \C$, let $\log z$ be the value of the logarithm of $z$ such that 
$0\le \arg z<2\pi$. We substitute 
$e^{(\swt w_{(3)}'-\swt w_{(1)}-\swt w_{(2)})\log z}$ for 
$x^{\swt w_{(3)}'-\swt w_{(1)}-\swt w_{(2)}}$
and  $(\log z)^k$ for $(\log x)^k$ for $k=0, \dots, k_{1}+k_{2}+k_{3}-3$
in $\langle w'_{(3)}, \mathcal{Y}(w_{(1)}, x)w_{(2)}\rangle$ to obtain a
complex number $\langle w'_{(3)}, \mathcal{Y}(w_{(1)}, z)w_{(2)}\rangle$.
Then we define $\langle w'_{(3)}, \mathcal{Y}(w_{(1)}, z)w_{(2)}\rangle$
for nonhomogeneous $w_{(1)}\in W_{1}$ and $w_{(2)}\in W_{2}$
and $w_{(3)}'\in W_{3}'$ using linearity. 

The $L(0)$-conjugation property  (\ref{l-0-conj-formal}) 
for logarithmic intertwining operators
recalled above gives
\begin{equation}\label{l-0-conj}
\langle w'_{(3)}, \mathcal{Y}(w_{(1)}, x)w_{(2)}\rangle
=\langle x^{L'(0)}w'_{(3)}, \mathcal{Y}(x^{-L(0)}w_{(1)}, 1)x^{-L(0)}w_{(2)}\rangle
\end{equation}
for $w_{(1)}\in W_{1}$, $w_{(2)}\in W_{2}$ and $w_{(3)}'\in W_{3}'$.

\begin{prop}\label{injectivity}
Assume 
that $W_{2}$ and $W_{3}'$ are generated by 
$\Omega^{0}_{N}(W_{2})$ and $\Omega^{0}_{N}(W_{3}')$, respectively.
Then the map 
\begin{eqnarray*}\rho: \mathcal{V}^{W_3}_{W_1\,W_2} &&\rightarrow\ \  {\rm Hom}_{A_N(V)}
(A_N(W_1)\otimes_{A_{N}(V)} \Omega^{0}_N(W_2), \Omega^{0}_N(W_3))\\
\mathcal {Y} &&\mapsto \ \ \rho(\mathcal {Y})
\end{eqnarray*}
is injective.
\end{prop}
\pf 
Assume that $\rho(\mathcal {Y}) = 0$. We have 
$W_3'=\coprod_{\mu\in\C/\Z}(W_{3}')^{\mu}$. For $w_{(2)} \in
(\Omega^{0}_N(W_2))_{[h_2]}$ and $w_{(3)}' \in (\Omega^{0}_N((W_3')^{\mu}))_{[h_3^{\mu}+n]}$,
\begin{eqnarray*}
\lefteqn{\langle w'_{(3)}, \mathcal {Y}(w_{(1)}, x)w_{(2)}\rangle}\nn
&&=\langle w'_{(3)}, (\pi^{\mu}\circ \mathcal {Y})(w_{(1)}, x)w_{(2)}\rangle\nn
&&=\langle x^{L'(0)}w'_{(3)}, (\pi^{\mu}\circ \mathcal {Y})(x^{-L(0)}w_{(1)},
1)x^{-L(0)}w_{(2)}\rangle\nn
&&=\langle x^{L'(0)}w'_{(3)}, \sum_{i\in \C}
(\pi^{\mu}\circ \mathcal {Y})_{i, 0}(x^{-L(0)}w_{(1)})x^{-L(0)}w_{(2)}\rangle\nn
&&=\langle x^{L'(0)}w'_{(3)}, 
(\pi^{\mu}\circ \mathcal {Y})_{\swt w_{(1)} + h_2 - h_3^{\mu}-n 
- 1, 0}(x^{-L(0)}w_{(1)})x^{-L(0)}w_{(2)}\rangle\nn
&&=\langle x^{L'(0)}w'_{(3)}, \sum_{\mu\in \C/\Z}
\sum_{m=0}^{N}\res_{x_{0}}x_{0}^{h_{2}-
h_3^{\mu}-m-1}
(\pi^{\mu}\circ \mathcal {Y})^{0}(x_{0}^{L(0)}x^{-L(0)}w_{(1)}, x_{0})
x^{-L(0)}w_{(2)}\rangle \nn
&&=\langle x^{L'(0)}w'_{(3)}, 
\sum_{\mu\in \C/\Z}\rho(\pi^{\mu}\circ \mathcal {Y})((x^{-L(0)}w_{(1)} 
+ O_N(W_1)) \otimes x^{-L(0)}w_{(2)})\rangle\nn
&&=\langle x^{L'(0)}w'_{(3)}, \rho(\mathcal {Y})((x^{-L(0)}w_{(1)} 
+ O_N(W_1)) \otimes x^{-L(0)}w_{(2)})\rangle\nn
&&= 0.
\end{eqnarray*}
Since $h_{2}$, $\mu$ and $n$ are arbitrary, this equality holds for 
$w_{(2)} \in
\Omega^{0}_N(W_2)$ and $w_{(3)}' \in \Omega^{0}_N(W_3')$.
Then by Proposition \ref{p3}, we have 
\[
\langle \tilde{w}'_{(3)}, \mathcal
{Y}(\tilde{w}_{(1)}, x)\tilde{w}_{(2)}\rangle = 0
\]
for all $\tilde{w}_{(1)} \in
W_1$, $\tilde{w}_{(2)} \in W_2$ and $\tilde{w}_{(3)}' \in
W_3'$. Thus $\mathcal {Y} = 0$. So $\rho$ is
injective.\epfv

\section{The main theorem}

In this section, under conditions stronger than 
those results in the preceding section, we state and 
prove our main theorem. 
The conditions needed in our main theorem is that
some lower-bounded generalized $V$-modules involved should satisfy
a certain universal property. Before we state and prove our main theorem,
we first give a construction of such lower-bounded generalized 
$V$-modules.

Let $W$ be a lower-bounded generalized $V$-module such that 
$W=\coprod_{n \in h_{W}+\N} W_{[n]}$ for some $h_W\in \C$
and $W_{[h_{W}]}\ne 0$.
Then  $G_{N}(W)=W_{[h_{W}+N]}$ 
is an $A_{N}(V)$-submodule 
of $\Omega^{0}_{N}(W)$. We now would like 
to construct a generalized $V$-module from $G_{N}(W)$ satisfying 
a certain universal property. We consider the affinization 
$V[t, t^{-1}]=V\otimes \C[t, t^{-1}]$ of $V$.
For simplicity, we shall use $u(m)$ to denote $u
\otimes t^m$
for $u \in V$ and $m \in \mathbb{Z}$. We consider the tensor algebra 
$T(V[t, t^{-1}])$. For
simplicity we shall omit the tensor product symbol when we write elements of 
$T(V[t, t^{-1}])$. 

For any $u \in V, m \in \mathbb{Z}$, $u(m)$ acts from the left on
$T(V[t, t^{-1}]) \otimes G_{N}(W)$. We shall also omit the tensor 
product symbol when we write elements of 
$T(V[t, t^{-1}])\otimes G_{N}(W)$.  The gradings on $T(V[t, t^{-1}])$
and $G_{N}(W)$ give a grading on 
$T(V[t, t^{-1}]) \otimes G_{N}(W)$. Explicitly, 
for homogeneous $u_i \in V, m_i
\in \mathbb{Z}, i= 1, \dots, s$ and homogeneous
$w\in G_{N}(W)$, the weight of 
$$u_1(m_1)\cdots
u_s(m_s)w$$
is
$${\rm wt}\ u_1 - m_1 - 1 + \cdots + {\rm wt}\ u_{s} - m_{s} - 1
 + {\rm wt}\ w.$$

Let
$$Y_{T(V[t, t^{-1}]) \otimes G_{N}(W)}(u, x): 
T(V[t, t^{-1}]) \otimes G_{N}(W) \rightarrow (T(V[t, t^{-1}])
\otimes
G_{N}(W))[[x, x^{-1}]]$$
be defined by
$$Y_{T(V[t, t^{-1}]) \otimes G_{N}(W)}(u, x) 
=\sum_{m \in \mathbb{Z}} u(m)x^{- m - 1}$$
for  $u \in V$.

Let $I_{V; W}$ be the $T(V[t, t^{-1}])$-submodule of $T(V[t, t^{-1}])
\otimes G_N(W)$ generated by the elements 
$$u(\wt u-1)w-o(u)w$$
for homogeneous $u\in V$ and $w \in G_{N}(W)$, the elements
$$u_1(m_1)\cdots u_s(m_s)w$$
for
homogeneous $u_i \in V$, $m_i \ge {\rm wt}\ u_i - 1$ satisfying $\sum_{i=1}^{s}
m_{i}> \sum_{i=1}^{s}(\wt u_{i}-1)+N$, $w \in G_{N}(W)$
and the coefficients in $x_1$ and $x_2$
of 
\begin{eqnarray*}
\lefteqn{Y_{T(V[t, t^{-1}]) \otimes G_{N}(W)}(u, x_1)
Y_{T(V[t, t^{-1}]) \otimes G_{N}(W)}(v, x_2)w}\nn
&&- Y_{T(V[t, t^{-1}]) \otimes G_{N}(W)}(v, x_2)
Y_{T(V[t, t^{-1}]) \otimes G_{N}(W)}(u,
x_1)w\nn
&& - {\rm Res}_{x_0}x_2^{-1}\delta\left(\frac{x_1
-x_0}{x_2}\right)Y_{T(V[t, t^{-1}]) \otimes G_{N}(W)}(Y_{V}(u, x_0)v, x_2)w
\end{eqnarray*}
for $u, v \in V$ and $w \in T(V[t, t^{-1}]) \otimes
G_{N}(W).$

Let 
$$\tilde{S}(V;  W) = (T(V[t, t^{-1}]) \otimes G_{N}(W))/I_{V; 
W_2}.$$
We shall use
elements of $T(V[t, t^{-1}])\otimes G_{N}(W)$ to represent elements
of $\tilde{S}(V;  W)$. But note that these elements now satisfy
the following relations: 
$$u(\wt u-1)w=o(u)w$$
for homogeneous $u\in V$ and $w \in G_{N}(W)$, 
$$u_1(m_1)\cdots u_s(m_s)w=0$$
for homogeneous $u_i \in V$, $m_i \ge {\rm wt}\ u_i - 1$ satisfying $\sum_{i=1}^{s}
m_{i}> \sum_{i=1}^{s}(\wt u_{i}-1)+N$, $w \in G_{N}(W)$
and 
\begin{eqnarray*}
\lefteqn{Y_{T(V[t, t^{-1}]) \otimes G_{N}(W)}(u, x_1)
Y_{T(V[t, t^{-1}]) \otimes G_{N}(W)}(v, x_2)w}\nn
&&- Y_{T(V[t, t^{-1}]) \otimes G_{N}(W)}(v, x_2)
Y_{T(V[t, t^{-1}]) \otimes G_{N}(W)}(u,
x_1)w\nn
&& = {\rm Res}_{x_0}x_2^{-1}\delta\left(\frac{x_1
-x_0}{x_2}\right)Y_{T(V[t, t^{-1}]) \otimes G_{N}(W)}(Y_{V}(u, x_0)v, x_2)w
\end{eqnarray*}
for $u, v \in V$ and $w \in T(V[t, t^{-1}]) \otimes
G_{N}(W).$

The map $Y_{T(V[t, t^{-1}])\otimes G_{N}(W)}$ induces a map for
$\tilde{S}(V; W)$ and we shall denote it by $Y_{\tilde{S}(V; W)}$.
By the definition of $\tilde{S}(V; W)$, the commutator formula
for $Y_{\tilde{S}(V;  W)}$ holds. 
Using this commutator formula
and other properties given by the definition of $\tilde{S}(V; W)$,
we see that 
$\tilde{S}(V; W)$ is 
spanned by elements of the form
\begin{equation}\label{span-0}
u_1(m_1)\cdots u_s(m_s)w 
\end{equation}
for
homogeneous $u_i \in V$, $m_i\in\Z$,
$i = 1, \dots,
s$ satisfying 
\begin{eqnarray*}
&\wt u_{1}-m_{1}-1\ge \cdots \ge \wt u_{s}-m_{s}-1,&\\
&\wt u_{1}-m_{1}-1+\cdots +\wt u_{s}-m_{s}-1\ge -N,&
\end{eqnarray*}
$w \in G_{N}(W)$. The grading of $T(V[t, t^{-1}])\otimes 
G_{N}(W)$
induces a grading on $\tilde{S}(V;  W)$ such that the weight of the element
(\ref{span-0}) is 
$$\wt u_{1}-m_{1}-1+\cdots +\wt u_{s}-m_{s}-1+\wt w.$$
Thus the real parts of the weights of the elements of $\tilde{S}(V;  W)$
are bigger than or equal to $\Re{(h_{W})}$. In particular, for $u \in V$ and $w \in \tilde{S}(V;  W)$, 
$u(m)w = 0$ when $m$ is sufficiently large. 

Recall the $\Z$-graded Lie algebra
$\hat{V}$ of operators on $V$ of the form $(Y_{V})_{n}(u)$ for $u\in V$ and 
$n\in \Z$, equipped with the Lie bracket for operators in Section 3.
Since the commutator formula
for $Y_{\tilde{S}(V; W)}$ holds, $\tilde{S}(V; W)$ is a graded $U(\hat{V})$-module. 
Let $P_{N}(\hat{V})
=\coprod_{k>N}^{\infty}\hat{V}_{(- k)}\oplus \hat{V}_{(0)}$. Then $P_{N}(\hat{V})$
is a subalgebra of $\hat{V}$. We know that 
$G_{N}(W)$ is an $A_{N}(V)$-module and is therefore  a graded
module for $\hat{V}_{(0)}$.
Let $\hat{V}_{(- k)}$ for $k>N$ act on $G_{N}(W)$ trivially. Then 
$G_{N}(W)$ is a $P_{N}(\hat{V})$-module.
Let 
$U(\cdot)$ be the universal enveloping 
algebra functor from the category of Lie algebras to
the category of associative algebras.
Then $\tilde{S}(V; W)$ as a graded $U(\hat{V})$-module is a quotient of 
the graded $U(\hat{V})$-module
$$U(\hat{V})\otimes_{U(P_{N}(\hat{V}))}G_{N}(W).$$

Let $J_{V; W}$ be the graded $U(\hat{V})$-submodule of $\tilde{S}(V;
W)$ generated by the coefficients in $x$ 
$$Y_{\tilde{S}(V; W)}(L(-1)u, x)w - \frac{d}{dx}Y_{\tilde{S}(V; W)}(u, x)w$$
and the coefficients in $x_0$, $x_2$ of
\begin{eqnarray*}
&{\displaystyle Y_{\tilde{S}(V; W)}(Y_{V}(u, x_0)v, x_2)w - 
Y_{\tilde{S}(V; W)}(u, x_0+x_{2})Y_{\tilde{S}(V; W)}(v, x_2)w}& \\
&\quad \quad\quad \quad\quad \quad\quad \quad + {\displaystyle  
{\rm Res}_{x_1}x_0^{-1}\delta\left(\frac{x_2  - x_1}{-x_0}\right)Y_{\tilde{S}(V; W)}(v,
x_2)Y_{\tilde{S}(V; W)}(u, x_1)w,}&
\end{eqnarray*}
for $u, v \in V$ and $w \in \tilde{S}(V;  W)$.
Let 
$$S_{N}(G_{N}(W))= \tilde{S}(V; W)/J_{V;  W}.$$
Then $S_{N}(G_{N}(W))$ is also a graded $U(\hat{V})$-module. We shall 
sometimes still use
elements of $T(V[t, t^{-1}])\otimes \Omega_N(W)$ to represent elements
of $S_{N}(V; W)$. But note that these elements now satisfy
more relations than the elements of $\tilde{S}(V;  W)$ written in 
the same form.
 The vertex operator map $Y_{\tilde{S}(V; W)}$
induces a vertex operator map 
$$Y_{S_{N}(G_{N}(W))}: V\otimes S_{N}(G_{N}(W))\to 
(S_{N}(G_{N}(W)))[[x, x^{-1}]].$$

In general, for a lower-bounded generalized $V$-module $W$, we know that 
there exists $h_{\mu}\in \C$ for $\mu\in \C/\Z$ such that 
$W=\oplus_{\mu\in \C/\Z}W^{\mu}$ where $W^{\mu}$ for $\mu\in \C/\Z$
are grading-restricted
generalized $V$-modules such that 
$W^{\mu}=\coprod_{n\in h_{\mu}+\N}(W^{\mu})_{[n]}$.
For $\mu\in \C/\Z$, we have defined 
$G_{N}(W^{\mu})$, $S_{N}(G_{N}(W^{\mu}))$ and 
$Y_{S_{N}(G_{N}(W^{\mu}))}$ above.
Let 
\begin{eqnarray*}
G_{N}(W)&=&\coprod_{\mu\in \C/\Z}G_{N}(W^{\mu}),\\
S_{N}(G_{N}(W))&=&\coprod_{\mu\in \C/\Z}S_{N}(G_{N}(W^{\mu}))
\end{eqnarray*}
and let 
$$Y_{S_{N}(G_{N}(W))}: V\otimes S_{N}(G_{N}(W))\to 
(S_{N}(G_{N}(W)))[[x, x^{-1}]]$$
be the map given by 
$$Y_{S_{N}(G_{N}(W))}(u, x)(w^{\mu_{1}}+\cdots +w^{\mu_{m}})
=Y_{S_{N}(G_{N}(W^{\mu_{1}}))}(u, x)w^{\mu_{1}}+\cdots +
Y_{S_{N}(G_{N}(W^{\mu_{m}}))}(u, x)w^{\mu_{m}}$$
for $u\in V$, $\mu_{1}, \dots, \mu_{m}\in \C/\Z$
and $w^{\mu_{1}}\in W^{\mu_{1}}, \dots, 
w^{\mu_{m}}\in W^{\mu_{m}}$.

\begin{thm}
Let $W$ be a lower-bounded generalized $V$-module. 
The graded space $S_{N}(G_{N}(W))$ equipped with 
vertex operator map $Y_{S_{N}(G_{N}(W))}$
is a lower-bounded generalized $V$-module
such that 
$G_{N}(S_{N}(G_{N}(W)))$ is isomorphic to $G_{N}(W)$.
The lower-bounded generalized $V$-module $S_{N}(G_{N}(W))$ 
satisfies the following universal 
property: For any lower-bounded 
generalized $V$-module $\widetilde{W}$  and any $A_{N}(V)$-module
map $\phi: G_{N}(W)\to G_{N}(\widetilde{W})$, there is a unique 
homomorphism $\bar{\phi}: S_{N}(G_{N}(W))\to \widetilde{W}$ of generalized 
$V$-modules
such that $\bar{\phi}|_{G_{N}(W)}=
\phi$.
\end{thm}
\pf
We prove the theorem only in the case that $W=\coprod_{n\in h_{W}+\N}W_{[n]}$;
the general case follows.
By definition, the commutator formula and the associator formula for the vertex operator 
map $Y_{S_{N}(G_{N}(W))}$ holds. Thus the Jacobi identity holds. 
The other properties are clearly satisfied. 

Since $\tilde{S}(V; W)$ is spanned by elements of form (\ref{span-0}), 
we see that the subspace
of $\tilde{S}(V; W)$ spanned by homogeneous elements
of weight  $h_{W}+N$ is isomorphic to $G_{N}(W)$. Since the
grading on $S_{N}(G_{N}(W))$ is induced from the 
grading on $\tilde{S}(V; W)$, we see that 
$G_{N}(S_{N}(G_{N}(W)))$ is isomorphic to $G_{N}(W)$.

The universal property follows from the construction.
\epfv

\begin{rema}
{\rm In \cite{DLM}, given an $A_{N}(V)$-module $U$,
an $\N$-gradable weak $V$-module  
$\bar{M}_{N}(U)$
is constructed. In the case that $U=G_{N}(W)$, it can be shown 
using the universal properties for $\bar{M}_{N}(U)$ and for $S_{N}(G_{N}(W))$
that the generalized $V$-module  $S_{N}(G_{N}(W))$
constructed above is isomorphic to $\bar{M}_{N}(U)$.}
\end{rema}

We need the following definition:

\begin{defn}
{\rm For $k\in \Z_{+}$, we say that {\it the 
$L(0)$-block sizes of a generalized $V$-module $W$
are less than $k$}  if for any homogeneous 
$w \in W$, $(L(0) - \wt w)^{k}w = 0$. We say that the $L(0)$-block sizes
of a generalized $V$-module is bounded if  there exists $k \in
\mathbb{Z}_{+}$ such that the 
$L(0)$-block sizes of $W$ are less than or equal to $k$.}
\end{defn}

\begin{rema}\label{r1}
{\rm Let $W_1$, $W_2$ and $W_3$ be lower-bounded generalized 
$V$-modules such that $W_{1}=\coprod_{n \in h_{1}+\Z} (W_{1})_{[n]}$,
$W_{2}=\coprod_{n \in h_{2}+\Z} (W_{2})_{[n]}$ and 
$W_{3}=\coprod_{n \in h_{3}+\Z} (W_{3})_{[n]}$ for some $h_1, h_2, h_3\in \C$
and such that the  $L(0)$-block sizes of $W_{1}$, $W_{2}$ 
and $W_{3}$ are less than $k_1$, $k_2$ and $k_3$, respectively.
Then by 
Proposition \ref{log:logwt}, we have
\begin{eqnarray*}\mathcal{Y}(w_{(1)}, x)w_{(2)} \in
x^{h_3-h_1-h_2}W_3[[x, x^{-1}]]\oplus x^{h_3-h_1-h_2}W_3[[x,
x^{-1}]]\log
x\oplus\nn \cdots\oplus x^{h_3-h_1-h_2}W_3[[x, x^{-1}]](\log x)^{k_1+k_2+k_3-3}.\nn
\end{eqnarray*}}
\end{rema}

Let $W_{1}$ and $W_{2}$ be lower-bounded generalized $V$-modules.
Assume that there exist $h_{1}, h_{2}\in \C$
such that  $W_{1}=\coprod_{n \in h_{1}+\Z} (W_{1})_{[n]}$ and
$W_{2}=\coprod_{n \in h_{2}+\Z} (W_{2})_{[n]}$ and that 
there exist 
$k_{1}, k_{2}\in \Z_{+}$ such that the $L(0)$-block sizes 
of  $W_1$ and $W_2$ 
are less than $k_{1}$ and $k_{2}$, respectively. Let $h_{3}\in \C$ and 
$k_{3}\in \Z_{+}$.
Similarly 
to the construction of $S_{N}(G_{N}(W))$ above, we now construct from $W_{1}$ and 
$G_{N}(W_{2})$
a lower-bounded
generalized $V$-module $S_{N}(V; W_1, W_2)$ and a logarithmic intertwining 
operator $\mathcal {Y}_{t}$ of type 
${S_{N}(V; W_1, W_2)\choose W_1 S_{N}(G_{N}(W_{2}))}$  such that 
$S_{N}(V; W_1, W_2)=\coprod_{n \in h_{3}+\Z} (S_{N}(V; W_1, W_2))_{[n]}$
and the $L(0)$-block sizes of $S_{N}(V; W_1, W_2)$ are less than or equal to 
$k_{3}$. 

Let $h = h_3 - h_1 - h_2$ and $k_0 = k_1 + k_2 + k_3$.  From Remark \ref{r1}, 
if $S_{N}(V; W_1, W_2)$ and $\mathcal {Y}_{t}$ are constructed,  we must have 
\begin{eqnarray*}
&\mathcal{Y}_{t}(w_{(1)}, x)w_{(2)} \in x^h S_{N}(V; W_1, W_2)
[[x, x^{-1}]]\oplus x^h S_{N}(V; W_1, W_2)[[x,
x^{-1}]]\log
x&\\
&\quad\quad\quad\quad
\oplus \cdots\oplus x^h S_{N}(V; W_1, W_2)[[x, x^{-1}]](\log x)^{k_0 - 3}.&
\end{eqnarray*}

We consider 
\[
t^{-h}W_1[t, t^{-1}][\log t] = W_1 \otimes t^{-h} \mathbb{C}[t, t^{-1}][\log t].
\]
For simplicity, we shall use $u(m)$ and $w_{(1)}(n, k)$ to denote 
$$u \otimes t^m\in V[t, t^{-1}]$$ and 
$$w_{(1)} \otimes t^n \otimes (\log t)^k\in
t^{-h}W_1[t, t^{-1}][\log t],$$
respectively,
for $u \in V, w_{(1)} \in W_1, m \in \mathbb{Z}, n \in -h + \mathbb{Z}$
and $k \in \mathbb{N}$. We consider the tensor algebra 
$$T(V[t,
t^{-1}]\oplus t^{-h}W_1[t, t^{-1}][\log t]).$$ 
The tensor algebra
$T(V[t, t^{-1}])$ is a subalgebra of this tensor algebra
and $t^{-h}W_1[t, t^{-1}][\log t]$
is a subspace. Let $T_{V; W_1}$ be the $T(V[t,
t^{-1}])$-sub-bimodule of $T(V[t, t^{-1}]\oplus t^{-h}W_1[t,
t^{-1}][\log t])$ generated by $t^{-h}W_1[t, t^{-1}][\log t]$.  Then
$T_{V; W_1} \otimes G_{N}(W_2)$ as a $T(V[t, t^{-1}])$-module is equivalent to 
$$T(V[t, t^{-1}])\otimes t^{-h}W_1[t, t^{-1}][\log t] \otimes T(V[t, t^{-1}])
\otimes G_{N}(W_2).$$
For
simplicity we shall omit the tensor product symbol for elements.
In particular, $T_{V; W_1} \otimes G_{N}(W_2)$ is spanned by elements of the form
$$u_1(m_1)\cdots u_s(m_s)w_{(1)}(n, k)u_{s+1}(m_{s+1})\cdots
u_{s+t}(m_{s+t})w_{(2)},$$
for $u_i \in V$, $m_i \in \mathbb{Z}$, $k \in
\mathbb{N}$, $i = 1, \dots, s+t$, $w_{(1)} \in W_1$ and $w_{(2)} \in
G_{N}(W_2)$.

For any $u \in V, m \in \mathbb{Z}$, $u(m)$ acts from the left on
$T_{V; W_1} \otimes G_{N}(W_2)$. The natural grading on $t^{-h}W_1[t, t^{-1}]$
gives a grading on $t^{-h}W_1[t, t^{-1}][\log t]$ with the weight of $\log t$ 
being $0$. The gradings on $T(V[t, t^{-1}])$,
$t^{-h}W_1[t, t^{-1}][\log t]$ and $G_{N}(W_2)$ give a grading on 
$T_{V; W_1} \otimes G_{N}(W_2)$. Explicitly, 
for homogeneous $u_i \in V, m_i
\in \mathbb{Z}, i= 1, \dots, s+t$, and homogeneous $w_{(1)} \in W_1$ and
 $w_{(2)} \in G_{N}(W_2)$, the weight of 
$$u_1(m_1)\cdots
u_s(m_s)w_{(1)}(n, k)u_{s+1}(m_{s+1})\cdots u_{s+t}(m_{s+t})w_{(2)}$$
is
$${\rm wt}\ u_1 - m_1 - 1 + \cdots + {\rm wt}\ u_{s+t} - m_{s+t} - 1
+ {\rm wt}\ w_{(1)} - n - 1 + {\rm wt}\ w_{(2)}.$$

For  $u \in V$, recall the map
$$Y_{T(V[t, t^{-1}]) \otimes G_N(W_2)}(u, x): 
T(V[t, t^{-1}]) \otimes G_N(W_2) \rightarrow (T(V[t, t^{-1}])
\otimes
G_N(W_2))[[x, x^{-1}]].$$
For $u\in V$  and $w_{(1)} \in W_1$, let
\begin{eqnarray*}
Y_{T_{V; W_1} \otimes G_N(W_2)}(u, x):& 
T_{V; W_1} \otimes G_N(W_2) &\rightarrow\; (T_{V;
W_1} \otimes
G_N(W_2))[[x, x^{-1}]],\nn
\mathcal {Y}_t(w_{(1)}, x):& T(V[t, t^{-1}]) \otimes G_{N}(W_2)
&\rightarrow\; x^h(T_{V; W_1} \otimes G_{N}(W_2))[[x,
x^{-1}]][\log x]
\end{eqnarray*}
be defined by
\begin{eqnarray*}
Y_{T_{V; W_1} \otimes G_N(W_2)}(u, x) &=&\sum_{m \in \mathbb{Z}} u(m)x^{- m - 1},\nn
\mathcal {Y}_t(w_{(1)}, x) &=& \sum_{n \in -h + \mathbb{Z}}\sum_{k
\in \mathbb{N}} w_{(1)}(n, k)x^{- n - 1}(\log x)^k,
\end{eqnarray*}
respectively.  

Let $I_{V; W_1, W_2}$ be the $T(V[t, t^{-1}])$-submodule of $T_{V;
W_1} \otimes G_N(W_2)$ generated by elements of the following
forms:
\begin{description}

\item $a u(\wt u-1)w_{(2)}-a o(u)
w_{(2)}$ for $a\in T_{V; W_{1}}$, $u \in V$, 
and $w_{(2)} \in G_{N}(W_{2})$,
 
\item $a u_1(m_1)\cdots u_s(m_s)w_{(2)}$
for  $a\in T_{V; W_{1}}$,
homogeneous $u_i \in V$, $m_i \ge {\rm wt}\ u_i - 1$ satisfying $\sum_{i=1}^{s}(\wt u_{i}
-m_{i}-1)<-N$, $w_{(2)} \in G_{N}(W_{2})$,

\item $u_1(m_1)\cdots u_s(m_s)w_{(1)}(n, k)w_{(2)}$  for homogeneous $u_i \in V$, 
$m_i\in \Z$, 
$w_{(1)} \in W_1$,  $n\in \C$, $k \in \N$, 
$w_{(2)} \in G_{N}(W_2)$, and either
$n \not\in -h + \mathbb{Z}$ or $\sum_{i=1}^{s}
\wt u_{i}-m_{i}-1+{\rm wt}\ w_{(1)} - n -1 + {\rm wt}\ w_{(2)} < \Re(h_3)$
or $k >k_{0}-3$,

\end{description}
and the coefficients in $x_1$, $x_2$ and $\log x_2$
of 
\begin{eqnarray*} 
\lefteqn{aY_{T(V[t, t^{-1}]) \otimes G_N(W_2)}(u, x_1)
Y_{T(V[t, t^{-1}]) \otimes G_N(W_2)}(v, x_2)w}\nn
&&- aY_{T(V[t, t^{-1}]) \otimes G_N(W_2)}(v, x_2)
Y_{T(V[t, t^{-1}]) \otimes G_N(W_2)}(u,
x_1)w\nn
&&- {\rm Res}_{x_0}x_2^{-1}\delta\left(\frac{x_1
-x_0}{x_2}\right)aY_{T(V[t, t^{-1}]) \otimes G_N(W_2)}(Y_{V}(u, x_0)v, x_2)w, \\
\lefteqn{Y_{T_{V; W_1} \otimes G_N(W_2)}(u, x_1)
Y_{T_{V; W_1} \otimes G_N(W_2)}(v, x_2)\tilde{w}}\nn
&&- Y_{T_{V; W_1} \otimes G_N(W_2)}(v, x_2)Y_{T_{V; W_1} \otimes G_N(W_2)}(u,
x_1)\tilde{w}\nn
&&- {\rm Res}_{x_0}x_2^{-1}\delta\left(\frac{x_1
-x_0}{x_2}\right)Y_{T_{V; W_1} \otimes G_N(W_2)}(Y_{V}(u, x_0)v, x_2)\tilde{w},
 \\
\lefteqn{Y_{T_{V; W_1} \otimes G_N(W_2)}(u, x_1)\mathcal {Y}_t(w_{(1)}, x_2)w}\nn
&&- \mathcal {Y}_t(w_{(1)},
x_2)Y_{T(V[t, t^{-1}]) \otimes G_N(W_2)}(u, x_1)w\nn
&&- {\rm Res}_{x_0}x_2^{-1}\delta\left(\frac{x_1
-x_0}{x_2}\right)\mathcal {Y}_t(Y_{W_{1}}(u, x_0)w_{(1)}, x_2)w,
\end{eqnarray*}
for $a\in T_{V; W_{1}}$, $u, v \in V, w_{(1)} \in W_1$,
$w\in T(V[t, t^{-1}] \otimes G_N(W_2)$ and $\tilde{w} \in T_{V;
W_1} \otimes G_N(W_2)$. 

Let 
$$\tilde{S}(V; W_1, W_2) = (T_{V; W_1} \otimes G_{N}(W_2))/I_{V; W_1,
W_2}.$$
We shall use
elements of $T_{V; W_1}\otimes G_{N}(W_2)$ to represent elements
of $\tilde{S}(V; W_1, W_2)$. But these elements now satisfy
relations. Recall the $U(\hat{V})$-module $\tilde{S}(V;  W_2)$ above.
The maps $Y_{T_{V; W_1} \otimes G_N(W_2)}(u, x)$ for 
$u\in V$ and $\mathcal{Y}_{t}(w_{(1)}, x)$ for 
$w_{(1)}\in W_{1}$ induce maps from 
$\tilde{S}(V; W_1, W_2)$ to $(\tilde{S}(V; W_1, W_2))[[x, x^{-1}]]$ and 
from $\tilde{S}(V;  W_2)$ to $x^{h}(\tilde{S}(V; W_1, W_2))[[x, x^{-1}]][\log x]$, 
respectively.
We shall use the notation $Y_{\tilde{S}(V; W_1, W_2)}(u, x)$ to denote 
the first map and the same notation $\mathcal{Y}_{t}(w_{(1)}, x)$ to denote the 
second map. These maps for all $u\in V$ and $w_{(1)}\in W_{1}$ give us  maps 
$Y_{\tilde{S}(V; W_1, W_2)}$ and $\mathcal{Y}_{t}$. Recall the map
$Y_{\tilde{S}(V;  W_2)}(u, x)$ for  $u \in V$ and $Y_{\tilde{S}(V;  W_2)}$.
By the definition of $\tilde{S}(V; W_1, W_2)$,  the commutator formulas
for $Y_{\tilde{S}(V; W_1, W_2)}$ and 
for $Y_{\tilde{S}(V; W_1, W_2)}$,  $\mathcal{Y}_{t}$ and $Y_{\tilde{S}(V;  W_2)}$
hold. Using these commutator formulas
and other properties given by the definition of $\tilde{S}(V; W_1, W_2)$,
we see that 
$\tilde{S}(V; W_1, W_2)$ is 
spanned by elements of the form
\begin{equation}\label{span}
u_1(m_1)\cdots u_s(m_s)w_{(1)}(n, k)w_{(2)} 
\end{equation}
for
homogeneous $u_i \in V$, $m_i\in\Z$, 
$i = 1, \dots,
s$, homogeneous $w_{(1)} \in W_1$, 
$n\in -h+\N$,  $0\le k\le k_{0}-3$ and 
$w_{(2)} \in G_{N}(W_2)$, satisfying 
\begin{eqnarray}
\Re(h_{3})&\le&  {\rm wt}\ w_{(1)} -n - 1 + {\rm wt}\ 
w_{(2)},\label{c2}\\
\Re(h_{3})&\le& \sum_{i=1}^{s}(\wt u_{i}-m_{i}-1)+\wt w_{(1)}-n-1+\wt w_{(2)}.
\label{c3}
\end{eqnarray}
The grading on $T_{V; W_1} \otimes G_N(W_2)$
induces a grading on $\tilde{S}(V; W_1, W_2)$ such that the weight of the element
(\ref{span}) is 
$$\wt u_{1}-m_{1}-1+\cdots +\wt u_{s}-m_{s}-1+\wt w_{(1)}-n-1+\wt w_{(2)}.$$
Thus the real parts of the weights of the elements of $\tilde{S}(V; W_1, W_2)$
are bigger than or equal to $\Re(h_{3})$.  In particular, for $u \in V, w_{(1)} \in
W_1$ and $w \in \tilde{S}(V; W_1, W_2)$, $u(m)w = 0$, $w_{(1)}(n, k)w =
0$ when $m$, $n$ and $k$ are sufficiently large. 

Since the commutator formulas
for $Y_{\tilde{S}(V; W_1, W_2)}$ holds, 
$\tilde{S}(V; W_1, W_2)$ is a $U(\hat{V})$-module. 

Let $J_{V; W_1, W_2}$ be the $U(\hat{V})$-submodule of $\tilde{S}(V;
W_1, W_2)$ generated by the coefficients in $x$ and $\log x$ of
\begin{eqnarray*}
&{\displaystyle aY_{\tilde{S}(V;  W_2)}(L(-1)u, x)w - a\frac{d}{dx}Y_{\tilde{S}(V;  W_2)}(u, x)w,}&\nn
&{\displaystyle Y_{\tilde{S}(V; W_1, W_2)}(L(-1)u, x)\tilde{w}
- \frac{d}{dx}Y_{\tilde{S}(V; W_1, W_2)}(u, x)\tilde{w},}&\nn
&{\displaystyle \mathcal {Y}_t(L(-1)w_{(1)}, x)w - \frac{d}{dx}\mathcal {Y}_t(w_{(1)},
x)w,}&\nn
&{\displaystyle \mathcal {Y}_t(w_{(1)}, x)w - x^{L(0)}\mathcal {Y}_t(x^{-L(0)}w_{(1)},
1)x^{-L(0)}w,}&\nn
\end{eqnarray*}
and the coefficients in $x_0$, $x_2$ and $\log x_2$ of
\begin{eqnarray*}
&{\displaystyle aY_{\tilde{S}(V;  W_2)}(Y_{V}(u, x_0)v, x_2)w 
- aY_{\tilde{S}(V;  W_2)}(u, x_0+x_{2})Y_{\tilde{S}(V;  W_2)}(v, x_2)w}& \\
&\quad \quad\quad \quad\quad \quad\quad \quad + {\displaystyle  
{\rm Res}_{x_1}x_0^{-1}\delta\left(\frac{x_2  - x_1}{-x_0}\right)aY_{\tilde{S}(V;  W_2)}(v,
x_2)Y_{\tilde{S}(V;  W_2)}(u, x_1)w,}&\\
&{\displaystyle Y_{\tilde{S}(V; W_1, W_2)}(Y_{V}(u, x_0)v, x_2)\tilde{w} - 
Y_{\tilde{S}(V; W_1, W_2)}(u, x_0+x_{2})Y_{\tilde{S}(V; W_1, W_2)}(v, x_2)\tilde{w}}& \\
&\quad \quad\quad \quad\quad \quad\quad \quad + {\displaystyle  
{\rm Res}_{x_1}x_0^{-1}\delta\left(\frac{x_2  - x_1}{-x_0}\right)
Y_{\tilde{S}(V; W_1, W_2)}(v,
x_2)Y_{\tilde{S}(V; W_1, W_2)}(u, x_1)\tilde{w},}&\\
&{\displaystyle \mathcal {Y}_t(Y_{V}(u, x_0)w_{(1)}, x_2)w - 
Y_{\tilde{S}(V; W_1, W_2)}(u, x_0+x_{2})\mathcal
{Y}_t(w_{(1)}, x_2)w}& \\
&{\displaystyle 
\quad \quad\quad \quad\quad \quad\quad \quad \quad
+ {\rm Res}_{x_1}x_0^{-1}\delta\left(\frac{x_2  - x_1}{-x_0}\right)\mathcal
{Y}_t(w_{(1)}, x_2)Y_{\tilde{S}(V;  W_2)}(u, x_1)w,}&
\end{eqnarray*}
for $a\in T_{V; W_{1}}$, $u, v \in V$, $w_{(1)} \in W_1$,
$w\in \tilde{S}(V;  W_2)$ and $\tilde{w} \in \tilde{S}(V; W_1, W_2)$.

Let 
$$S_{N}(V; W_1, W_2) = \tilde{S}(V; W_1, W_2)/J_{V; W_1, W_2}.$$
Then $S_{N}(V;
W_1, W_2)$ is also a $U(\hat{V})$-module. We shall use
elements of $T_{V; W_1}\otimes G_N(W_2)$ to represent elements
of $S_{N}(V; W_1, W_2)$. But  these elements now satisfy
more relations than the elements of $\tilde{S}(V; W_1, W_2)$ written in 
the same form.

The maps $Y_{\tilde{S}(V; W_1, W_2)}(u, x)$ for $u\in V$
and $\mathcal{Y}_{t}(w_{(1)}, x)$
for $w_{(1)}\in W_{1}$ induce maps 
$$Y_{S_{N}(V; W_1, W_2)}(u, x): S_{N}(V; W_1, W_2)\to 
(S_{N}(V; W_1, W_2))[[x, x^{-1}]]$$ 
and 
$$\mathcal {Y}_t(w_{(1)}, x): S_{N}(G_{N}(W_2)) \to 
x^{h}(S_{N}(V; W_1, W_2))[[x, x^{-1}]][\log x],$$
respectively. These maps give maps
\begin{eqnarray*}
Y_{S_{N}(V; W_1, W_2)}: V\otimes S_{N}(V; W_1, W_2)&\to& 
(S_{N}(V; W_1, W_2))[[x, x^{-1}]]\nn
u\otimes w&\mapsto&  Y_{S_{N}(V; W_1, W_2)}(u, x)w
\end{eqnarray*}
and 
\begin{eqnarray*}
\mathcal {Y}_t: W_{1}\otimes S_{N}(G_{N}(W_2)) &\to &
x^{h}(S_{N}(V; W_1, W_2))[[x, x^{-1}]][\log x]\nn
w_{(1)}\otimes w_{(2)}&\mapsto& \mathcal {Y}_t(w_{(1)}, x)w_{(2)}.
\end{eqnarray*}
By construction, these operators satisfy the lower truncation property,
the identity property for $Y_{S_{N}(V; W_1, W_2)}$, the commutator formula for 
$Y_{S_{N}(V; W_1, W_2)}$
and for $Y_{S_{N}(V; W_1, W_2)}$,
$\mathcal{Y}_t$ and $Y_{S_{N}(G_{N}(W_2))}$, the associator formula for
$Y_{S_{N}(G_{N}(W_2))}$ and for $Y_{S_{N}(V; W_1, W_2)}$,
$\mathcal{Y}_t$ and $Y_{S_{N}(G_{N}(W_2))}$, 
the $L(-1)$-derivative property for
$Y_{S_{N}(V; W_1, W_2)}$ and for $\mathcal{Y}_t$ 
and the $x^{L(0)}$-conjugation property for 
$Y_{S_{N}(V; W_1, W_2)}$ and for
$\mathcal {Y}_t$.

\begin{thm}\label{s-n-v-w-1-w-2}
Let $W_{1}$ and $W_{2}$ be lower-bounded generalized $V$-modules.
Assume that there exist $h_{1}, h_{2}\in \C$
such that  $W_{1}=\coprod_{n \in h_{1}+\Z} (W_{1})_{[n]}$ and
$W_{2}=\coprod_{n \in h_{2}+\Z} (W_{2})_{[n]}$ and that 
there exist 
$k_{1}, k_{2}\in \Z_{+}$ such that the $L(0)$-block sizes 
of  $W_1$ and $W_2$ 
are less than $k_{1}$ and $k_{2}$, respectively. Let $h_{3}\in \C$ and 
$k_{3}\in \Z_{+}$. Then 
the graded space $S_{N}(V; W_1, W_2)$ equipped with 
vertex operator map $Y_{S_{N}(V; W_1, W_2)}$
is a lower-bounded generalized $V$-module
such that $S_{N}(V; W_1, W_2)=\coprod_{n \in h_{3}+\Z} (S_{N}(V; W_1, W_2))_{[n]}$
and the $L(0)$-block sizes of $S_{N}(V; W_1, W_2)$ are less than or equal to 
$k_{3}$ and $\mathcal {Y}_{t}$ is a logarithmic intertwining 
operator of type 
${S_{N}(V; W_1, W_2)\choose W_1 S_{N}(G_{N}(W_{2}))}$. 
\end{thm}
\pf 
Since $Y_{S_{N}(V; W_1, W_2)}$ satisfies the lower truncation property,
the identity property, the commutator formula,
the associator formula and 
the $L(-1)$-derivative property, $S_{N}(V; W_1, W_2)$ is a 
generalized $V$-module. 
Since $\mathcal {Y}_t$ satisfies the lower truncation property,
the commutator formula,
the associator formula and 
the $L(-1)$-derivative property,
$\mathcal {Y}_{t}$ is a logarithmic intertwining 
operator of type 
${S_{N}(V; W_1, W_2)\choose W_1 S_{N}(G_{N}(W_{2}))}$.
It is clear from the construction and the properties above satisfied by 
$Y_{S_{N}(V; W_1, W_2)}$ that $S_{N}(V; W_1, W_2)$
is lower bounded, 
$S_{N}(V; W_1, W_2)=\coprod_{n \in h_{3}+\Z} (S_{N}(V; W_1, W_2))_{[n]}$
and the $L(0)$-block sizes of $S_{N}(V; W_1, W_2)$ are less than or equal to 
$k_{3}$. 
\epfv

Now we state and prove our main result:

\begin{thm}\label{main}
Let $W_1$, $W_2$ and  $W_3$ be lower-bounded generalized $V$-modules 
whose $L(0)$-block sizes are bounded.
Let $N$ be a nonnegative integer such that  $W_1$ is generated by
$\Omega^{0}_{N}(W_1)$ and $W_2$ and  $W'_3$ are isomorphic to 
$S_{N}(G_{N}(W_2))$ and $S_{N}(G_{N}(W'_3))$, respectively. 
Then the map
\begin{eqnarray*}
\rho: \mathcal{V}^{W_3}_{W_1\,W_2} &\rightarrow&
{\rm Hom}_{A_N(V)}
(A_N(W_1)\otimes_{A_N(V)}\Omega_{N}^{0}(W_2), \Omega_{N}^{0}(W_3))\\
\mathcal {Y} &\mapsto &\rho(\mathcal {Y})
\end{eqnarray*}
is a linear isomorphism.
\end{thm}
\pf 
Since $S_{N}(G_{N}(W_2))$ and $S_{N}(G_{N}(W'_3))$ are generated by 
$G_{N}(S_{N}(G_{N}(W_2)))\subset \Omega^{0}_{N}(S_{N}(G_{N}(W_2)))$ and 
$G_{N}(S_{N}(G_{N}(W'_3)))\subset \Omega^{0}_{N}(S_{N}(G_{N}(W'_3)))$, 
respectively, they are generated by $\Omega^{0}_{N}(S_{N}(G_{N}(W_2)))$ and 
$\Omega^{0}_{N}(S_{N}(G_{N}(W'_3)))$, respectively. Since
$W_{2}$ and $W'_{3}$ are isomorphic to  
$S_{N}(G_{N}(W_2))$ and $S_{N}(G_{N}(W'_3))$, respectively,
$W_{2}$ and $W'_{3}$ 
are generated by $\Omega^{0}_{N}(W_{2})$ and $\Omega^{0}_{N}(W'_{3})$, 
respectively.
By Proposition \ref{injectivity}, $\rho$ is injective. 
So we need only prove that $\rho$ is surjective. Given any
element $f$ of  
\[
{\rm Hom}_{A_N(V)}
(A_N(W_1)\otimes_{A_N(V)}\Omega_{N}^{0}(W_2), \Omega_{N}^{0}(W_3)),
\] 
we want to
construct an element $\mathcal {Y}^f$ of $\mathcal{V}^{W_3}_{W_1\,W_2}$
such that $\rho(\mathcal {Y}^f) = f$. 

We first construct $\mathcal {Y}^f$ in the case that there exist $h_{i}\in \C$
for $i=1, 2, 3$ such that $W_{1}=\coprod_{n \in h_{1}+\Z} (W_{1})_{[n]}$,
$W_{2}=\coprod_{n \in h_{2}+\Z} (W_{2})_{[n]}$ and 
$W_{3}=\coprod_{n \in h_{3}+\Z} (W_{3})_{[n]}$. Since the $L(0)$-block sizes 
of  $W_1$, $W_2$ and  $W_3$ are bounded, there exists positive integers
$k_{1}$, $k_{2}$ and $k_{3}$ such that the $L(0)$-block sizes 
of  $W_1$, $W_2$ and  $W_3$ 
are less than $k_{1}$, $k_{2}$ and $k_{3}$, respectively.
Let $h = h_3 - h_1 - h_2$ and $k_0 = k_1 + k_2 + k_3$. Then by
Proposition \ref{r1}, for the logarithmic intertwining operator $\mathcal {Y}^f$
that we want to construct, we have 
\begin{eqnarray*}
\mathcal{Y}(w_{(1)}, x)w_{(2)} \in x^h W_3[[x, x^{-1}]]\oplus x^h W_3[[x,
x^{-1}]]\log
x\oplus \cdots\oplus x^h W_3[[x, x^{-1}]](\log x)^{k_0 - 3}.
\end{eqnarray*}

From Theorem \ref{s-n-v-w-1-w-2}, we have a lower-bounded generalized $V$-module
$S_{N}(V; W_1, W_2)$
such that $S_{N}(V; W_1, W_2)=\coprod_{n \in h_{3}+\Z} (S_{N}(V; W_1, W_2))_{[n]}$
and the $L(0)$-block sizes of $S_{N}(V; W_1, W_2)$ are less than or equal to 
$k_{3}$ and a logarithmic intertwining 
operator $\mathcal {Y}_{t}$ of type 
${S_{N}(V; W_1, W_2)\choose W_1 S_{N}(G_{N}(W_{2}))}$. 
Let $S_{N}^{0}(V; W_1, W_2)$ be generalized $V$-submodule 
of $S_{N}(V; W_1, W_2)$ generated by elements 
of the form 
$w_{(1)}(n, 0)w_{(2)}\in S_{N}(V; W_1, W_2)$
where 
$w_{(1)} \in W_1$ is homogeneous, $n \leq {\rm wt}\ w_{(1)} - 1 + {\rm wt}\
w_{(2)} - h_3$ and $w_{(2)} \in G_{N}(W_2)$ is homogeneous.
Then  we have a homomorphism $\mu : S_{N}^{0}(V; W_1, W_2)
\rightarrow W_3$ of generalized $V$-modules 
defined as follows: 


Using commutator formulas, 
we know that $S_{N}^{0}(V; W_1, W_2)$
is spanned by elements of the form (\ref{span})
for
homogeneous $u_i \in V$, $m_i\in\Z$, 
$i = 1, \dots,
s$, homogeneous $w_{(1)} \in W_1$, 
$n\in -h+\N$,  and 
$w_{(2)} \in G_{N}(W_2)$, satisfying (\ref{c2}) and (\ref{c3}).
For $w'_{(3)}\in G_{N}(W'_{3})$, 
we define an element $\mu'(w_{(3)}')\in (S_{N}^{0}(V; W_1, W_2))'$
by
\begin{eqnarray}\label{mu'}
\lefteqn{\langle \mu'(w_{(3)}'), u_1(m_1)\cdots u_s(m_s)
w_{(1)}(n, 0)w_{(2)}\rangle}\nn
&&=\langle w_{(3)}', 
(Y_{W_{3}})_{m_1}(u_1)\cdots (Y_{W_{3}})_{m_s}(u_s)P_{\swt w_{(1)}+\swt w_{(2)}-n-1}
f((w_{(1)} + O_N(W_{1}))\otimes w_{(2)})\rangle.\nn
\end{eqnarray}
Note that the only relations among elements of form (\ref{span}) are
those given by $I_{V; W_{1}, W_{2}}$ and $J_{V; W_{1}, W_{2}}$.
These relations are also satisfied by elements of $W_{3}$
of the form 
$$(Y_{W_{3}})_{m_1}(u_1)\cdots (Y_{W_{3}})_{m_s}(u_s)P_{\swt w_{(1)}+\swt w_{(2)}-n-1}
f((w_{(1)} + O_N(W_{1}))\otimes w_{(2)}).$$
Hence $\mu'(w_{(3)}')$ is well defined. More precisely, we can see that 
$\mu'(w_{(3)}')$ is well defined as follows:
Consider the graded subspace of 
$T_{V;W_{1}}\otimes G_{N}(W_{2})$ spanned by elements
of the form 
$$u_1(m_1)\cdots u_s(m_s)
w_{(1)}(n, 0)w_{(2)}$$
for homogeneous $u_i \in V$, $m_i\in\Z$, 
$i = 1, \dots,
s$, homogeneous $w_{(1)} \in W_1$, 
$n\in -h+\N$,  and 
$w_{(2)} \in G_{N}(W_2)$, satisfying (\ref{c2}) and (\ref{c3}).
First,  we define  $\mu'(w_{(3)}')$ using (\ref{mu'}) to be 
an element of the graded dual space of this graded subspace of 
$T_{V;W_{1}}\otimes G_{N}(W_{2})$.
Then from the definitions of $\mu'(w_{(3)}')$ and 
$I_{V; W_{1}, W_{2}}$, we see that
$\mu'(w_{(3)}')$ annihilates the intersection of $I_{V; W_{1}, W_{2}}$
and this graded subspace of 
$T_{V;W_{1}}\otimes G_{N}(W_{2})$. So 
$\mu'(w_{(3)}')$ is in fact an element of the graded dual space 
of the graded subspace of $\tilde{S}(V; W_{1}, W_{2})$
spanned by elements of the form 
$$u_1(m_1)\cdots u_s(m_s)
w_{(1)}(n, 0)w_{(2)}$$ 
for homogeneous $u_i \in V$, $m_i\in\Z$, 
$i = 1, \dots,
s$, homogeneous $w_{(1)} \in W_1$, 
$n\in -h+\N$,  and 
$w_{(2)} \in G_{N}(W_2)$, satisfying (\ref{c2}) and (\ref{c3}).
But from the definitions of 
$\mu'(w_{(3)}')$ and $J_{V; W_{1}, W_{2}}$,  we see that
$\mu'(w_{(3)}')$ annihilates  the intersection of $J_{V; W_{1}, W_{2}}$ and this 
graded subspace of $\tilde{S}(V; W_{1}, W_{2})$. Thus 
$\mu'(w_{(3)}')$ is in fact an element of $(S_{N}^{0}(V; W_1, W_2))'$.

By definition, we see that if $w_{(3)}'$ is homogeneous, 
then $\mu'(w_{(3)}')$ is also homogeneous and $\wt \mu'(w_{(3)}')=\wt w_{(3)}'$.
Thus 
$\mu'(w_{(3)}')\in G_{N}((S_{N}^{0}(V; W_1, W_2))')$ for any $w_{(3)}'
\in G_{N}(W'_{3})$
and we obtain a linear map $\mu': G_{N}(W'_{3})\to 
G_{N}((S_{N}^{0}(V; W_1, W_2))')$. The map $\mu'$ is 
in fact a homomorphism of $A_{N}(V)$-modules, that is,
$$\mu'((Y_{W_{3}}')_{\swt u-1}(u)w'_{(3)})=
(Y_{(S_{N}^{0}(V; W_1, W_2))'})_{\swt u-1}(u)\mu'(w'_{(3)}),$$
for homogeneous $u\in V$ and $w'_{(3)}\in G_{N}(W'_{3})$. But this is 
equivalent to 
$$\mu'((Y_{W_{3}'}^{o})_{\swt u-1}(u)w'_{(3)})=
(Y_{(S_{N}^{0}(V; W_1, W_2))'}^{o})_{\swt u-1}(u)\mu'(w'_{(3)})$$
for homogeneous $u\in V$ and $w'_{(3)}\in G_{N}(W'_{3})$,
which follows from the calculation
\begin{eqnarray*}
\lefteqn{\langle \mu'((Y_{W_{3}'}^{o})_{\swt u-1}(u)w_{(3)}'), 
u_1(m_1)\cdots u_s(m_s)
w_{(1)}(n, 0)w_{(2)}\rangle}\nn
&&=\langle (Y_{W_{3}'}^{o})_{\swt u-1}(u)w_{(3)}', 
(Y_{W_{3}})_{m_1}(u_1)\cdots (Y_{W_{3}})_{m_s}(u_s)\cdot\nn
&&\quad\quad\quad\quad\quad\quad\quad\quad\quad\quad\quad\quad
 \cdot P_{\swt w_{(1)}+\swt w_{(2)}-n-1}
f((w_{(1)} + O_N(W_{1}))\otimes w_{(2)})\rangle\nn
&&=\langle w_{(3)}', (Y_{W_{3}})_{\swt u-1}(u)
(Y_{W_{3}})_{m_1}(u_1)\cdots (Y_{W_{3}})_{m_s}(u_s)\cdot\nn
&&\quad\quad\quad\quad\quad\quad\quad\quad\quad\quad\quad\quad
 \cdot P_{\swt w_{(1)}+\swt w_{(2)}-n-1}
f((w_{(1)} + O_N(W_{1}))\otimes w_{(2)})\rangle\nn
&&=\langle \mu'(w_{(3)}'), u(\wt u-1)
u_1(m_1)\cdots u_s(m_s)
w_{(1)}(n, 0)w_{(2)}\rangle\nn
&&=\langle (Y_{(S_{N}^{0}(V; W_1, W_2))'}^{o})_{\swt u-1}(u)\mu'(w_{(3)}'), 
u_1(m_1)\cdots u_s(m_s)
w_{(1)}(n, 0)w_{(2)}\rangle
\end{eqnarray*}
for
homogeneous $u_i \in V$, $m_i\in\Z$, 
$i = 1, \dots,
s$, homogeneous $w_{(1)} \in W_1$, 
$n\in -h+\N$,  and 
$w_{(2)} \in G_{N}(W_2)$, satisfying (\ref{c2}) and (\ref{c3}).

Since 
$W'_{3}$ is isomorphic to $S_{N}(G_{N}(W'_3))$,
by the universal property of $S_{N}(G_{N}(W'_3))$,
we obtain a homomorphism, still denoted as $\mu'$, of generalized $V$-modules from 
$W'_{3}$ to $(S_{N}^{0}(V; W_1, W_2))'$, extending $\mu': G_{N}(W'_{3})\to 
G_{N}((S_{N}^{0}(V; W_1, W_2))')$. The adjoint map $\mu''$ of $\mu'$ 
is a homomorphism of generalized $V$-modules from $(S_{N}^{0}(V; W_1, W_2))''$
to $W_{3}$. In particular, the restriction of $\mu''$ to $S_{N}^{0}(V; W_1, W_2)$
is a homomorphism $\mu$ 
of generalized $V$-modules from $S_{N}^{0}(V; W_1, W_2)$
to $W_{3}$.

We now define 
$$(\mathcal
{Y}^f)_{n, 0}(w_{(1)})w
=\mu(w_{(1)}(n, 0)w)$$
for
$w_{(1)} \in W_1$,  
$w\in S_{N}(G_{N}(W_{2}))$ and 
$n\in -h+\N$. Since $W_{2}$ is isomorphic to $S_{N}(G_{N}(W_{2}))$,
we obtain maps $(\mathcal
{Y}^f)_{n, 0}(w_{(1)}): W_{2} \to W_{3}$ for $n\in -h+\N$.
Let 
$$(\mathcal
{Y}^f)^{0}(w_{(1)}, x)=\sum_{n\in -h+\N}(\Y^f)_{n, 0}(w_{(1)})x^{-n-1}.$$
In particular, 
for $w_{(1)}\in W_{1}$, $w_{(2)}\in W_{2}$ and $w'_{(3)}\in W'_{3}$,
$$\langle w'_{(3)}, (\Y^{f})^{0}(w_{(1)}, 1)w_{(2)}\rangle
=\sum_{n\in -h+\N}\langle w'_{(3)},(\Y^f)_{n, 0}(w_{(1)})w_{(2)}\rangle$$
is well defined. 
Now we construct a logarithmic intertwining operator $\mathcal
{Y}^f$ of type ${W_3\choose W_1\,W_2}$ by
$$\langle w'_{(3)}, \Y^{f}(w_{(1)}, x)
w_{(2)}\rangle=
\langle x^{L'(0)}w'_{(3)}, (\Y^{f})^{0}(x^{-L(0)}w_{(1)}, 1)
x^{-L(0)}w_{(2)}\rangle$$
for $w_{(1)}\in W_{1}$, $w_{(2)}\in W_{2}$ and $w'_{(3)}\in W'_{3}$ (recall
(\ref{l-0-conj})).

Since $\mathcal {Y}_t$
satisfies the commutator formula, the associative formula and the
$L(-1)$-derivative property, so does $\mathcal {Y}^f$. Thus
$\mathcal {Y}^f$ satisfies the Jacobi identity and the
$L(-1)$-derivative property. So it is a logarithmic intertwining
operator of the desired type. It is clear from the construction that
$\rho(\mathcal {Y}^f) = f$. 

In the general case, by Remark \ref{congruent},
$W_{1}$, $W_{2}$ and $W_{3}'$ can all be 
decomposed as direct sums of grading generalized modules $W$ such that 
$G_{N}(W)$ 
are spanned by homogeneous elements of weights $h+n$ for $n=0, \dots, N$.
Then the logarithmic intertwining
operator $\mathcal {Y}^f$ can be obtained by adding those intertwining 
operators obtained from the case discussed above.
\epfv

From Proposition \ref{p1} and Theorem \ref{main}, we obtain immediately the following
result:

\begin{cor}
Assume that  there exists a positive integer 
$N$ such that the absolute value of the difference of the lowest weights of
any two irreducible $V$-modules is less than or equal to $N$. 
Let $W_1$, $W_2$ and  $W_3$ be grading-restricted generalized $V$-modules of finite lengths
whose $L(0)$-block sizes are bounded.
Assume that $W_2$ and  $W'_3$ are equivalent to 
$S_{N}(G_{N}(W_2))$ and $S_{N}(G_{N}(W'_3))$, respectively. 
Then the map
\begin{eqnarray*}
\rho: \mathcal{V}^{W_3}_{W_1\,W_2} &\rightarrow&
{\rm Hom}_{A_N(V)}
(A_N(W_1)\otimes_{A_N(V)}\Omega_{N}^{0}(W_2), \Omega_{N}^{0}(W_3))\\
\mathcal {Y} &\mapsto &\rho(\mathcal {Y})
\end{eqnarray*}
is a linear isomorphism.\epf
\end{cor}

\noindent {\small \sc Department of Mathematics, Rutgers University,
110 Frelinghuysen Rd., Piscataway, NJ 08854-8019}

\noindent {\it and}

\noindent {\small \sc Beijing International Center for Mathematical Research,
Peking University, Beijing, China}
\vspace{1em}

\noindent {\em E-mail address}: yzhuang@math.rutgers.edu, 

\vspace{1em}

\noindent {\small \sc Department of Mathematics, Rutgers University,
110 Frelinghuysen Rd., Piscataway, NJ 08854-8019}

\vspace{1em}

\noindent {\em E-mail address}: yookinwi@math.rutgers.edu

\end{document}